\documentclass{aims}

\usepackage{txfonts}
\usepackage{amsmath}
\usepackage{amssymb}
\usepackage{bbm}
\usepackage{algorithm}
\usepackage{comment}
\usepackage{t1enc}
\usepackage{bm}
\usepackage{xr}
\usepackage[noend]{algpseudocode}

\numberwithin{equation}{section}

\DeclareMathOperator*{\argmin}{arg\,min}

\def\typeofarticle{Research article}
\def\currentvolume{x}
\def\currentissue{x}
\def\currentyear{2022}

\def\ppages{xxx--xxx}
\def\DOI{}
\def\Received{January 2022}
\def\Accepted{September 2022}
\def\Published{December 2022}

\begin{document}
	
\title{Comparison of the performance and reliability between improved sampling strategies for polynomial chaos expansion}

\author{%
	Konstantin Weise\affil{1,2,+}\corrauth,
	Erik M\"uller\affil{1,2,+},
	Lucas Po\ss{}ner\affil{1,3}
	and
	Thomas R. Kn\"osche\affil{1}	
}

\shortauthors{the Author(s)}

\address{%
	\addr{\affilnum{1}}{Max Planck Institute for Human Cognitive and Brain Sciences, Brain Networks Group, Stephanstra\ss{}e 1a, 04103, Leipzig, Germany}
	\addr{\affilnum{2}}{Technische Universit\"at Ilmenau, Advanced Electromagnetics Research Group, Helmholtzplatz 2, 98693 Ilmenau, Germany}
	\addr{\affilnum{3}}{Hochschule f\"ur Technik Wirtschaft und Kultur Leipzig, W\"achterstra\ss{}e 13, 04107 Leipzig, Germany}
	\addr{\affilnum{+}}{Contributed equally}}

\corraddr{konstantin.weise@tu-ilmenau.de; Tel: +49-341-9940-2580; Fax:\\ +49-341-9940-2624.
}

\begin{abstract}
As uncertainty and sensitivity analysis of complex models grows ever more important, the difficulty of their timely realizations highlights a need for more efficient numerical operations. Non-intrusive Polynomial Chaos methods are highly efficient and accurate methods of mapping input-output relationships to investigate complex models. There is substantial potential to increase the efficacy of the method regarding the selected sampling scheme. We examine state-of-the-art sampling schemes categorized in space-filling-optimal designs such as Latin Hypercube sampling and L1-optimal sampling and compare their empirical performance against standard random sampling. The analysis was performed in the context of L1 minimization using the least-angle regression algorithm to fit the GPCE regression models. Due to the random nature of the sampling schemes, we compared different sampling approaches using statistical stability measures and evaluated the success rates to construct a surrogate model with relative errors of $<0.1$\%, $<1$\%, and $<10$\%, respectively. The sampling schemes are thoroughly investigated by evaluating the y of surrogate models constructed for various distinct test cases, which represent different problem classes covering low, medium and high dimensional problems. Finally, the sampling schemes are tested on an application example to estimate the sensitivity of the self-impedance of a probe that is used to measure the impedance of biological tissues at different frequencies.  We observed strong differences in the convergence properties of the methods between the analyzed test functions.

\end{abstract}


\maketitle
	
\section{Introduction}\label{Intro}

\subsection{Background}
During the last several decades, computational advances as well as new insights that reduce the computational cost of methods have made many problems feasible to solve. This has a great impact for uncertainty and sensitivity analysis as a means of reliability engineering, where a huge number of evaluations need to be carried out to thoroughly investigate their behavior and to identify critical model parameters. In recent years, substantial development has advanced the current landscape of stochastic computational problems. Surrogate modeling has become a very popular method to effectively map complex input-output relationships and to tremendously reduce the computational cost of successive uncertainty and sensitivity evaluations. The general Polynomial Chaos Expansion (GPCE) \cite{Wiener.1938, Xiu.2002, LeMaitre.2010} has proven to be a versatile method to predict model behavior for systems where direct evaluation can be cumbersome and very time consuming \cite{Eldred.2012}. GPCE allows physical and engineering systems to be investigated with various Quantities Of Interest (QOI) as functions of various uncertain model parameters. The goal lies in characterizing the sensitivities of the output dependencies of the system inputs, which can be done by performing a Sobol decomposition \cite{Sobol.1999} or gradient based measures \cite{Xiu.2009}. A lot of research can be found on this process of uncertainty quantification (UQ) using the GPCE model \cite{Ghanem.2003, Xiu.2002, LeMaitre.2010}. It has found applications in non-destructive material testing \cite{Weise.2016}, neuroscience \cite{Codecasa.2016, Weise.2015, Saturnino.2019b, Weise.2020b}, mechanical engineering \cite{Sepahvand.2012}, aerospace engineering\cite{Hosder.2007}, electrical engineering \cite{Kaintura.2018, Diaz.2018}, fluid dynamics \cite{Najm.2009}, and various other fields. However, compared to engineering, there is a relative lack of targeted applications of GPCE in the life sciences \cite{Burk.2020, Grignard.2022, Hu.2018, Massoud.2019, Pepper.2019, Son.2020}. In this branch of science, strong assumptions are made about the parameter values to be chosen. Moreover, the analysis of individuals requires the analysis of model behavior to consider stochastic parameter definitions instead of deterministic approaches.  

\subsection{Recent advancements}
The majority of modern research on GPCE has focused on non-intrusive approaches where the problem at hand can be treated as a black box system. There exist a number of possibilities to further improve the basic GPCE approach. On one hand, it is possible to modify the assembly process of the basis functions by identifying and choosing the most suitable order \cite{Sachdeva.2006}, splitting the GPCE problem in a multi-element-GPCE (ME-GPCE) \cite{Wan.2006}, or applying an adaptive algorithm to extend the number of samples and basis functions iteratively \cite{Blatman.2011, Novak.2021}. This reduces the number of unknown GPCE coefficients and the number of required model evaluations.

\subsection{Improved sampling strategies}
On the other hand, the potential of a more efficient GPCE approximation lies in the selection of the sampling locations prior to any GPCE approximation. To take a closer look at this topic, we thoroughly investigate GPCE optimized sampling schemes compared to standard random sampling. We will (i) show the performance of standard Monte Carlo methods in the framework of polynomial chaos, (ii) improve them with space-filling sampling designs using state-of-the-art Latin Hypercube Sampling (LHS) schemes \cite{McKay.1979, Helton.2003}, (iii) apply principles of Compressed Sensing (CS) and optimal sampling for L1-minimization \cite{Donoho.2006, Rauhut.2012, Hampton.2015b}. We investigate the performance and reliability of the sampling schemes in a comprehensive numerical study, which consists of three representative test cases and one practical example.
On one hand, LHS has seen a lot of usage in the fields of reliability engineering and uncertainty analysis since it is equipped with good space-filling properties \cite{Shu.2011, Robinson.1999}. The basic principle was improved by Jin et al. (2003) \cite{Jin.2003}, who proposed an optimization scheme based on the Enhanced Stochastic Evolutionary algorithm (ESE), which maximizes the Maximum-Minimal Distance criterion \cite{Johnson.1990} to reliably construct sample sets with a very even spread. It can be reasonably assumed that GPCE can benefit significantly from this, since it ensures to some extent that the parameter space is scanned evenly and thus all features of the transfer function can be found.
On the other hand, CS emerged in the field of efficient data recovery to reconstruct signals with a much smaller number of samples than the Shannon-Nyquist criterion would suggest \cite{Donoho.2006, Candes.2008, Candes.2006, Elad.2007}. This has been applied in a number of cases where the number of samples available is limited \cite{Lustig.2005, Paredes.2007, Ender.2010}. Because it is not possible to select the required basic functions in advance, most GPCE dictionaries are over-complete, which leads to sparse coefficient vectors. Using these properties, compressive sampling recently became popular in the framework of UQ. Another appealing fact is that in computational UQ, it is possible to freely draw additional samples, thus enabling a multitude of new possibilities compared to real data acquisition, where the number of measurements is limited and possibly even restricted. This provided a new subcategory of sparse Polynomial Chaos Expansions \cite{Blatman.2011, Karagiannis.2014, Jakeman.2015, Deman.2016} where generally solvers like Least Angle Regression (LARS) \cite{Tibshirani.2004} or Orthogonal Matching Pursuit (OMP) \cite{Pati.1993} are used to determine the sparse coefficient vectors. Better GPCE recoverability has also been shown by designing a unique sampling method for Legendre Polynomials using the Chebyshev measure \cite{Rauhut.2012} or by extending this with a distinct coherence parameter and sampling from an altered input variable distribution by using a Monte-Carlo Markov-Chain algorithm \cite{Hampton.2015b}. Those methods, however, are restricted to problems with a low number of random variables employing a high polynomial order. Progress has also been made in defining criteria such as the mutual coherence $\mu$, the RIP-constant \cite{Candes.2008b}, and a number of correlation-constants to categorize measurement matrices and quantify their possible recovery. In this paper, we focus on evaluating the mutual coherence parameter as a global measure for a minimization objective and on a combination of different local criteria. We adopted the proposed "near optimal" sampling method of Alemazkoor and Meidani \cite{Alemazkoor.2018}, which uses a greedy algorithm to ensure a more stable recovery. We additionally use the same framework to create mutual coherence optimal GPCE matrices. This is meant to serve as a comparative example. Additionally, we propose a global approach to create an L1-optimal design by using an iterative algorithm to maximize the local and global optimality criteria.
Given those two approaches to construct the set of sampling points, we are proposing a hybrid design that is partially created using LHS and then expanded according to a chosen L1-optimality criteria or vice versa. This aims to give a broad overview over not only the effectiveness of the two branches, but also over possible enhancements or coupling occurrences from their interaction.
In this paper, we are going to compare the aforementioned sampling strategies on a set of test problems with varying order and dimension. We investigate their error convergence over different sample sizes. We also test their applicability on a practical example, which consists of an electrode-impedance model used to characterize the impedance of brain tissue. All algorithms are implemented in the open source python package "pygpc" \cite{Weise.2020}, and the scripts to run the presented benchmarks are provided in the Supplemental Material.
The remainder of the paper is structured as follows. 

\subsection{Content}
The theoretical background of GPCE is revisited in Section \ref{sec:PolynomialChaosExpansion}. It is followed by introducing the different sampling schemes, namely standard Random Sampling in Section \ref{sec:RandomSampling}, LHS in Section \ref{sec:Space-fillingoptimalsamplingapproaches} and CS-optimal sampling in Section \ref{sec:CompressiveSamplingapproaches}. An overview about the test problems to which the sampling schemes are applied is given in Section \ref{sec:Results} together with the benchmark results \ref{sec:Results}. Finally, the results are discussed in Section \ref{sec:Discussion}.

\section{Polynomial Chaos Expansion}
\label{sec:PolynomialChaosExpansion}
In GPCE, the $d$ parameters of interest, which are assumed to underlie a distinct level of uncertainty, are modeled as a $d$-variate random vector denoted by $\bm{\xi} = (\xi_1, \, \xi_2, \, ... \xi_d)$ following some probability density function (pdf) $p_i(\xi_i)$, with $i=1,...,d$. The random parameters are defined in the probability space $(\Theta, \Sigma, P)$. The event or random space $\Theta$ contains all possible events. $\Sigma$ is a $\sigma$-Algebra over $\Theta$, containing sets of events, and $P$ is a function assigning the probabilities of occurrence to the events. The number of random variables $d$ determines the dimension of the uncertainty problem. It is assumed that the parameters are statistically mutually independent from each other. In order to perform a GPCE expansion, the random variables must have a finite variance, which defines the problem in the $L_2$-Hilbert space. 

The quantity of interest (QOI) will be analyzed in terms of the random variables $\bm{\xi}$, is $y(\mathbf{r})$. It may depend on some external parameters $\mathbf{r}=(r_{0},\,...,\,r_{R-1})$ like space, where $R=3$, or any other dependent parameters. Those are treated as deterministic and are not considered in the uncertainty analysis. 

The basic concept of GPCE is to find a functional dependence between the random variables $\bm{\xi}$ and the solutions $y(\mathbf{r},\bm{\xi})$ by means of an orthogonal polynomial basis $\Psi(\bm{\xi})$. In its general form, it is given by: 
\begin{equation}
y(\mathbf{r},\bm{\xi}) = \sum_{\bm{\alpha}\in\mathcal{A}} c_{\bm{\alpha}}(\mathbf{r}) \Psi_{\bm{\alpha}}(\bm{\xi}).
\label{eq:ua:gPC_coeff_form}
\end{equation}

A separate GPCE expansion must be performed for every considered parameter set $\mathbf{r}$. The discrete number of QOIs is denoted as $N_y$.

The terms are indexed by the multi-index $\bm{\alpha}=(\alpha_0,...,\alpha_{d-1})$, which is a $d$-tuple of non-negative integers $\bm{\alpha}\in\mathbb{N}_0^d$. The sum is carried out over the multi-indices, contained in the set $\mathcal{A}$. 

The function $\Psi_{\bm{\alpha}}(\bm{\xi})$ are the polynomial basis functions of GPCE. They are composed of polynomials $\psi_{\alpha_i}(\xi_i)$.
\begin{equation}
\Psi_{\bm{\alpha}}(\bm{\xi}) = \prod_{i=1}^{d} \psi_{\alpha_i}(\xi_i)
\label{eq:ua:Psi}
\end{equation}

The polynomials $\psi_{\alpha_i}(\xi_i)$ are defined for each random variable separately according to the corresponding input pdf $p_i(\xi_i)$. They must be chosen to be orthogonal with respect to the pdfs of the random variables, e.g. Jacobin polynomials for beta-distributed random parameters or Hermite polynomials for normal-distributed random variables. 
The family of polynomials for an optimal basis of continuous probability distributions is given by the Askey scheme \cite{Askey.1985}. The index of the polynomials denotes its order (or degree). In this way, the multi-index $\bm{\alpha}$ corresponds to the order of the individual basis functions forming the joint basis function.

In general, the set $\mathcal{A}$ of multi-indices can be freely chosen according to the problem under investigation. In practical applications, the \emph{maximum order} GPCE is frequently used. In this case, the set $\mathcal{A}$ includes all polynomials whose total order does not exceed a predefined order $p$. 
In the present work, the concept of \emph{maximum order} GPCE is extended by introducing the \emph{interaction} order $p_i$. An interaction order $p_i(\bm{\alpha})$ can be assigned to each multi-index $\bm{\alpha}$. The multi-index reflects the respective powers of the polynomial basis functions of random variables, $\Psi_{\bm{\alpha}}(\bm{\xi})$:

\begin{align}
p_i(\bm{\alpha}) = \lVert\bm{\alpha}\rVert_0,
\end{align}

where $\lVert\bm{\alpha}\rVert_0 = \#(i:\alpha_i>0)$ is the zero (semi)-norm, quantifying the number of non-zero index entries. The reduced set of multi-indices is then constructed by the following rule:

\begin{equation}
\mathcal{A}(p, p_i) := \left\{ \bm{\alpha} \in \mathbb{N}_0^d\, : \lVert \bm{\alpha} \rVert_1 \leq p \wedge \lVert\bm{\alpha}\rVert_0 \leq p_i \right\} 
\label{eq:ua:multi_index_set_reduced}
\end{equation}

It includes all elements from a total order GPCE with the restriction of the interaction order $p_i$. Reducing the number of basis functions is advantageous especially in case of high-dimensional problems. This is supported by observations in a number of studies, where the magnitude of the coefficients decreases with increasing order and interaction \cite{Hampton.2015}. Besides that, no hyperbolic truncation was applied to the basis functions \cite{Blatman.2011}.

After constructing the polynomial basis, the corresponding GPCE-coefficients $c_{\bm{\alpha}}(\mathbf{r})$ must be determined for each output quantity. In this regard, the output variables are projected from the $d$-dimensional probability space $\Theta$ into the $N_c$-dimensional polynomial space $\mathcal{P}_{N_c}$. This way, an analytical approximation of the solutions $y(\mathbf{r},\bm{\xi})$ as a function of its random input parameters $\bm{\xi}$ is derived and very computationally-efficient investigation of its stochastics is made possible. 

The GPCE from (\ref{eq:ua:gPC_coeff_form}) can be written in matrix form as:

\begin{align}\label{eq:gPC_system}
\mathbf{Y} = \bm{\Psi}\mathbf{C}
\end{align}

Depending on the sampling strategy, one may define a diagonal positive-definite matrix $\mathbf{W}$ whose diagonal elements $\mathbf{W}_{i,i}$ are given by a function of sampling points $w(\bm{\xi}^{(i)})$.

\begin{align}\label{eq:gPC_system_weighted}
\mathbf{W}\mathbf{Y} = \mathbf{W}\bm{\Psi}\mathbf{C}
\end{align}

The GPCE-coefficients for each QOI (columns of $\mathbf{C}$) can then be found by using solvers that minimize either the L1 or the L2 norm of the residuum depending on the expected sparsity of the coefficient vectors. Each row in \ref{eq:gPC_system_weighted} corresponds to a distinct sampling point $\bm{\xi}_i$. For this reason, the choice of the sampling points has a considerable influence on the characteristics and solvability of the equation system.

%
%
%

Complex numerical models can be very computationally intensive. To enable uncertainty and sensitivity analysis of such models, the number of sampling points must be reduced to a minimum. Minimizing the sampling points may lead to a situation where there are fewer observations than unknowns, i.e. $M \leq K$, resulting in a under-determined system of equations with infinitely many solutions for $\mathbf{c}$. Considering compressive sampling, we want $\mathbf{c}$ to be the sparsest solution, formulating the recovery problem as:

\begin{equation}\label{eq:l0}
\min_{\mathbf{c}} ||\mathbf{c}||_0 \quad \textrm{subject to} \quad \mathbf{\Psi} \mathbf{c} = \mathbf{u}
\end{equation}

where $||.||_0$ indicates the $\ell_0$-norm, the number of non-zero entries in $\mathbf{c}$. This optimization problem is NP-hard and not convex. The latter property can be overcome by reformulating it using the L1-norm:

\begin{equation}\label{eq:l1_gpc}
\min_{\mathbf{c}} ||\mathbf{c}||_1 \quad \textrm{subject to} \quad \mathbf{\Psi} \mathbf{c} = \mathbf{u}
\end{equation}

It has been shown that if $[\mathbf{\Psi}]$ is sufficiently incoherent and $\mathbf{c}$ is sufficiently sparse, the solution of the $\ell_0$ minimization is unique and equal to the solution of the L1 minimization \cite{Bruckstein.2009}. The minimization in equation \ref{eq:l1_gpc} is called basis pursuit \cite{Chen.2001} and can be solved using linear programming.

\section{Sampling techniques}
\subsection{Standard Monte Carlo sampling}
\label{sec:RandomSampling}
The most straightforward sampling method is to draw samples according to the input distributions. In this case, one proceeds with a Monte Carlo method to sample the random domain without any sophisticated process for choosing the sampling locations. The random samples must be chosen independently and should be uncorrelated, but a simple sampling process may inadvertently violate this requirement, especially when the number of sampling points is small (a situation we are targeting). For instance, the sampling points can be concentrated in certain regions that do not reveal some important features of the model's behavior, thus significantly degrading the overall quality of the GPCE approximation.

\subsection{Coherence-optimal sampling}
\label{sec:COSampling}
Coherence-optimal (CO) sampling aims to improve the stability of the coefficients when solving (\ref{eq:gPC_system_weighted}). It was introduced by Hampton and Doostan in the framework of GPCE in \cite{Hampton.2015}. The Gram matrix (also referred to as the gramian or information matrix) defined in eq. (\ref{eq:gramian}) and its properties play a central role when determining the GPCE coefficients. Coherence-optimal sampling has been the building block for a number of sampling strategies that aim for an efficient sparse recovery of the PC \cite{Alemazkoor.2018, Hadigol.2018}. It generally outperforms random sampling by a large margin on problems with higher order than dimensionality $p \geq d $ and has been claimed to perform well on any given problem when incorporated in compressive sampling approaches \cite{Hampton.2015b}. It is defined by:

\begin{equation}\label{eq:gramian}
\mathbf{G_\Psi} = \frac{1}{N_g}\mathbf{\Psi^T}\mathbf{\Psi}
\end{equation}

CO sampling seeks to minimize the spectral matrix norm between the Gram matrix and the identity matrix, i.e. $||\mathbf{G_\Psi}-\mathbf{I}||$, by minimizing the coherence parameter $\mu$:

\begin{equation}\label{eq:coherence}
\mu = \sup_{\bm{\xi}\in\Omega} \sum_{j=1}^P \left|w(\mathbf{\xi})\psi_j(\bm{\xi})\right|^2
\end{equation} 

This can be done by sampling the input parameters with an alternative distribution:

\begin{equation}\label{eq:alternative_distribution}
P_{\mathbf{Y}}(\bm{\xi}) := c^2 P(\bm{\xi}) B^2(\bm{\xi}),
\end{equation} 

where $c$ is a normalization constant, $P(\bm{\xi})$ is the joint probability density function of the original input distributions, and $B(\bm{\xi})$ is an upper bound of the PC basis:

\begin{equation}\label{eq:B2}
B(\bm{\xi}):= \sqrt{\sum_{j=1}^P|\psi_j(\bm{\xi})|^2}
\end{equation}

To avoid defining the normalization constant $c$, a Markov Chain Monte Carlo approach using a Metropolis-Hastings sampler \cite{Hastings.1970} is used to draw samples from $P_{\mathbf{Y}}(\bm{\xi})$ in (\ref{eq:alternative_distribution}). For the Mertopolis-Hastings sampler, it is necessary to define a sufficient candidate distribution. For a coherence optimal sampling according to (\ref{eq:coherence}), this is realized by a proposal distribution $g(\xi)$ \cite{Hampton.2015}. By sampling from a different distribution than $P(\xi)$, however, it is not possible to guarantee $\mathbf{\Psi}$ to be a matrix of orthonormal polynomials. Therefore $\mathbf{W}$ needs to be a diagonal positive-definite matrix of weight-functions $w(\xi)$. In practice, it is possible to compute $\mathbf{W}$ with:

\begin{equation}\label{eq:weighting_matrix}
 w_i(\xi) = \frac{1}{B_i(\xi)}
\end{equation}

A detailed description about the technique can be found in \cite{Hampton.2015}.
	
\subsection{Optimal design of experiment}
\label{sec:DSampling}
A judicious choice of sampling points $\{\bm{\xi}^{(i)}\}_{i}^{N_g}$ allows us to improve the properties of the Gramian without any prior knowledge about the model under investigation. The selection of an appropriate optimization criterion derived from $[\mathbf{G_\Psi}]$ and the identification of the corresponding optimal sampling locations is the core concept of optimal design of experiment (ODE). The most popular criterion is $D$-optimality, where the goal is to increase the information content from a given number of sampling points by minimizing the determinant of the inverse of the Gramian:

\begin{equation}\label{eq:det_Dopt}
\phi_D = \left|\mathbf{G_\Psi}^{-1}\right|^{1/N_c}, 
\end{equation}

$D$-optimal designs are focused on precise estimation of the coefficients. Besides $D$-optimal designs, there exist many other alphabetic optimal designs, such as $A$-, $E$-, $I$-, or $V$- optimal designs with different goals and criteria. A nice overview of these designs can be found in \cite{Pukelsheim.2006, Atkinson.2007}.

Hadigol and Doostan investigated the convergence behavior of $A$-, $D$- and $E$-optimal designs \cite{Hadigol.2018} in combination with coherence-optimal sampling in the framework of least squares GPCE. They found that those designs clearly outperform standard random sampling. Their analysis was restricted to cases where the number of sampling points is larger than the number of unknown coefficients ($N_g > N_c$). Based on the current state of knowledge, our analysis focuses on investigating the convergence properties of $D$-optimal and $D$-coherence-optimal designs in combination with L1 minimization where $N_g < N_c$.
	
\subsection{Space-filling optimal sampling}
\label{sec:Space-fillingoptimalsamplingapproaches}
In order to overcome the disadvantages of standard random sampling for low sample sizes, one may use sampling schemes that improve the coverage of the random space. Early work on this topic focused on pseudo-random sampling while optimizing distinct distance and correlation criteria between the sampling points. Designs optimizing the Maximum Minimal distance \cite{Johnson.1990, Morris.1995, Jin.2003} or Audze-Eglais Designs \cite{Audze.1977, Bates.2004} proved to be both more efficient and more reliable than standard random sampling schemes. Space-filling optimal sampling such as Latin Hypercube Sampling (LHS) is nowadays frequently being used in the framework of GPCE \cite{Choi.2004, Hosder.2007, Hadigol.2018}. In the following, we briefly introduce two prominent distance criteria we used in our space-filling optimal sampling approaches.

\subsubsection{Space-filling optimality criteria}
\paragraph*{Maximum-Minimal distance criterion:}
The maximum-minimal distance criterion is a space-filling optimality criterion. A design can be called maximum-minimum distance optimal if it maximizes the minimum inter-site distance \cite{Johnson.1990}:

\begin{equation}\label{eq:inter_cite_dist}
\min_{1 \leq i, j\leq n, i\neq j} d(\mathbf{x}_{i},\mathbf{x}_{j})\quad \textrm{subject to} \quad d(\mathbf{x}_{i},\mathbf{x}_{j}) = d_{ij} = \left( \sum_{k=1}^{m}|x_{ik}-x_{jk}|^t\right)^\frac{1}{t}
\end{equation}

where $d(\mathbf{x}_{i},\mathbf{x}_{j})$ is the distance between two sampling points $\mathbf{x}_{i}$ and $\mathbf{x}_{j}$, and $t = 1$ or $2$. A design optimized in its minimum inter-site distance is able to create well-distributed sampling points. For a low number of sampling points, however, the sampling points may be heavily biased towards the edges of the sampling space because the distance criterion pushes the sampling points relentlessly outwards and away from possible features close to the center of the sampling space \cite{Morris.1995}.

\paragraph*{The $\varphi_{p}$ criterion:}
To counteract the shortcomings of the plain inter-side distance in (\ref{eq:inter_cite_dist}), the equivalent $\varphi_{p}$ criterion has been proposed \cite{Morris.1995}. A $\varphi_{p}$-optimal design can be constructed by setting up a distance list $(d_{1},...,d_{s})$ obtained by sorting the inter-site distances $d_{ij}$ together with a corresponding index list $(J_{1}, ..., J_{s})$. The $d_{i}$ are distinct distance values, $d_{1}<d_{2}<...<d_{s}$ and $J_{s}$ are the corresponding indices of pairs of sites in the design separated by $d_{i}$. A design can then be called $\varphi_{p}$-optimal if it minimizes:

\begin{equation}\label{eq:phi_p}
\varphi_{p} = \left(\sum_{i=1}^{s}J_{i}d_{i}^{-p}\right)^\frac{1}{p}
\end{equation}

We empirically choose $p$ as 10 in the numerical construction of LHS designs.

\subsection{Limitations of the distance criterion for low sampling sizes}\label{subsec:phi-limit}

As the investigation of the following sampling schemes in a sparse reconstruction suggests, functions with a high number of variables occur very commonly. Since the goal of a sparse reconstruction is to reduce the sampling size, an important caveat to the optimization of criteria based around the distance $d_i, d_{i,j}$ in the two criteria maximum-minimal distance and the $\phi_p$ is that it shows a systematic bias that breaks the uniformity sought in the following optimization algorithm in section \ref{sec:ESE}. This problem has been identified recently by Vořechovský and Eliáš \cite{Vorechovsky.2020, Elias.2020} and becomes apparent in efficient optimization. They introduced a new distance criterion called periodic distance:

\begin{equation}\label{eq:periodic-d}
	\bar{d_{ij}} = \left( \sum_{k=1}^{m} \left(\min \left(|x_{ik}-x_{jk}|, 1 - |x_{ik}-x_{jk}| \right) \right)^t \right)^\frac{1}{t}
\end{equation}

With $\bar{d_{ij}}$, it is possible to calculate the periodic maximum-minimum distance and the periodic $\phi_p$ criterion as $  \min_{1 \leq i, j\leq n, i\neq j} \bar{d_{ij}} $ and $ \phi_p(\bar{d_{ij}})$, respectively, by using the periodic distance instead of the conventional euclidean distance metric. Further, for the $\phi_p$ criterion, it was shown by the same authors and Mašek that specifying the $p$-exponent based on investigations of the potential energy of the design can lead to an enhancement in the space-filling and projection properties as well as a decrease in its discrepancy. For successful application of LHS, they recommend $p = N_{var} + 1$, where $N_{var}$ is the number of variables in the given function \cite{Vorechovsky.2020b}.

\subsubsection{Standard Latin Hypercube Sampling (LHS)}
In LHS, the $d$-dimensional sampling domain is segmented into $n$ subregions corresponding to the $n$ sampling points to be drawn. LHS designs ensure that every subregion is sampled only once. This method can be mathematically expressed by creating a matrix of sampling points $\mathbf{\Pi}$ with:

\begin{equation}\label{eq:lhs}
\pi_{i, j} = \frac{p_{i, j} - u_{i,j}}{n}, 
\end{equation}

where $\mathbf{P}$ is a matrix of column-wise randomly permuted indices of its rows and $\mathbf{U}$ is a matrix of independent uniformly distributed random numbers $u \in [0, 1]$.


\subsubsection{$\varphi_{p}$-optimal Latin Hypercube Sampling}
The space-filling properties of LHS designs can be improved by optimizing the $\varphi_{p}$ criterion. A pseudo-optimal design can be determined by creating a pool of $n_i$ standard LHS designs and choosing the one with the best $\varphi_{p}$ criterion. As $n_i$ reaches infinity, the design will become space-filling optimal. In this study we used $n_i=100$ iterations, which was found to be an efficient trade-off between computational cost and $\varphi_{p}$-optimality.

\subsubsection{Enhanced Stochastic Evolutionary Algorithm LHS}\label{sec:ESE}

The Enhanced Stochastic Evolutionary Algorithm Latin Hypercube Sampling (LHS-ESE) is a very stable space-filling optimal algorithm designed by Jin et al. (2003) \cite{Jin.2003}. The resulting designs aim for a specified $\varphi_p$ parameter and achieve that by multiple element-wise exchanges of an initial LHD in an inner loop while storing their respective $\varphi_p$ - parameters in an outer loop. This process shows a far smaller variance in the space-filling criteria of the created sample sets.

\begin{figure*}[tbh]
	\centering
	\includegraphics[width=0.49\textwidth]{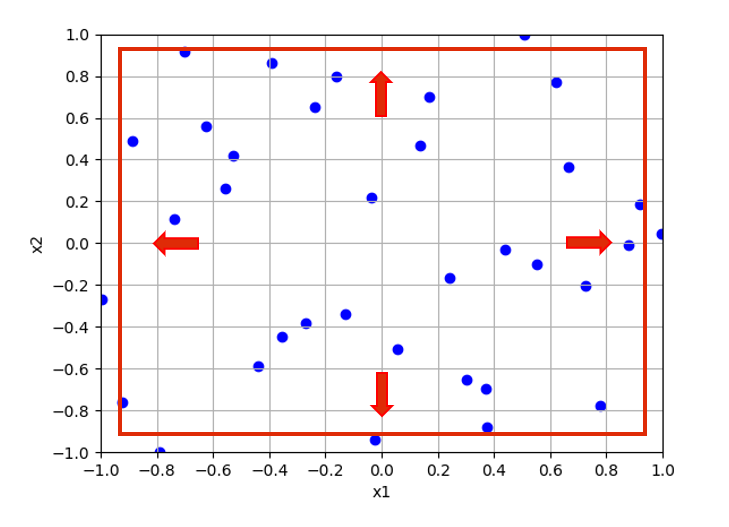}
	\caption[SC-ESE]{Schematic representation of the SC-ESE, where the outer area is cut to an $\alpha$ fraction of its original size and the center is stretched outward}
	\label{fig:SC_ESE_scheme}
\end{figure*}

However, we observed that the LHS-ESE scheme often undersamples the boundaries of the random domain, which is disadvantageous for transfer functions with high gradients close to the parameter boundaries. This is less apparent for a high number of samples but becomes a serious drawback when the number of sampling points is low; this effect can be linked to the systematic bias pointed out in \cite{Vorechovsky.2020}. In order to overcome this problem, we modified the LHS-ESE algorithm by shrinking the first and last subregion of the interval to a fraction of their original sizes while keeping the remaining intermediate $n-2$ subregions equally spaced. The procedure is illustrated in Fig. \ref{fig:SC_ESE_scheme}. The initial matrix for the Latin Hypercube design $\mathbf{\Pi}$ will then be changed as $\mathbf{P}$ is randomly permuted according to:

\begin{equation}\label{eq:sese_1}
\pi_{i, 1} = \pi_{i, n} =\alpha \frac{p_{i, j} - u_{i,j}}{n},
\end{equation}

where $\alpha$ is the fraction to which the size of the border interval is decreased, and $i \in [1, d]$ and $j \in [1, n]$ are used to cover the size reduction of the intervals at the edges. Our empirical studies showed that a reduction of $\alpha=\frac{1}{4}$ counteracts the aforementioned undersampling close to the border. We pick the index $j$ to target these edges, since after the normalization by dividing by $n$, the indices $1$ and $n$ for $j$ are expected for the values closest to $0$ and $1$ respectively, which occur at the borders of the sampled section. The centre is then stretched by:

\begin{equation}\label{eq:sese_2}
\scalebox{1.25}{$\pi_{i, j_c} = \begin{cases}
\frac{p_{i, j_c} - u_{i,j_c}}{n} - \frac{(1-\alpha) p_{i, j_c}}{n -2} &\text{for $j \leq \frac{n}{2}$}\\
\frac{p_{i, j_c} - u_{i,j_c}}{n} + \frac{(1-\alpha) p_{i, j_c}}{n -2} &\text{else,}
\end{cases}$}
\end{equation}

with $j_c$ being the indices of $j$ without the border domains $1$ and $n$, $j_c = j \setminus \{1, n\}$. After that alteration, the elements of each column in $\mathbf{\Pi}$ can be randomly permuted to proceed with the construction of the Latin Hypercube design just like in \ref{eq:lhs}. If $\alpha$ is made smaller, then the size of the guaranteed sampling region at the border also becomes smaller, thus forcing the sampling point to be chosen closer to the border as illustrated in Fig. \ref{fig:SC_ESE_scheme}. To the best of our knowledge, the Enhanced Stochastic Evolutionary LHS algorithm has not yet been studied in the context of GPCE.


\subsection{Compressive Sampling}
\label{sec:CompressiveSamplingapproaches}
Compressive Sampling is a novel method, first introduced in the field of signal processing, that allows the recovery of signals with significantly fewer samples assuming that the signals are sparse: i.e. a certain portion of the coefficients are zero, meaning that the coefficient vector $\mathbf{c}$ can be well-approximated with only a small number of non-vanishing terms. A coefficient vector $\mathbf{c}$ that is $s$-sparse obeys:

\begin{equation}\label{eq:sparsity:2}
||\mathbf{c}||_0 \leq s, \quad \forall s \in \mathbb{N}	
\end{equation}

The locations of the sampling points have a profound impact on the reconstruction quality because they determine the properties of the GPCE matrix. There are several criteria that can be evaluated exclusively on the basis of the GPCE matrix that may favor the reconstruction. It could be shown that optimization of those criteria lead to designs that promote successful reconstruction \cite{Gang.2013, Alemazkoor.2018}. In the following, we give a brief overview about the different criteria we considered in this study.

\subsubsection{L1-optimality criteria}
 	
\paragraph*{Mutual Coherence}
The mutual coherence (MC) of a matrix measures the cross-correlations between its columns by evaluating the largest absolute and normalized inner product between different columns. It can be evaluated by:

\begin{equation}\label{eq:L1:mc}
\mu(\mathbf{\Psi}) = \max_ {1 \leq i, j\leq N_c, j\neq i} \quad \frac{|\psi_i^T \psi_j|}{||\psi_i||_2||\psi_j||_2}
\end{equation} 

The objective is to select sampling points that minimize $\mu(\mathbf{\Psi})$ for a desired L1-optimal design. It is noted that minimizing the mutual-coherence considers only the worst-case scenario and does not necessarily improve compressive sampling performance in general \cite{Elad.2007}.

%
%

\paragraph*{Average Cross-Correlation}
It is shown in [39, 45, 46, 47] that the robustness and accuracy of signal recovery can be increased by minimizing the distance between the Gram matrix $\mathbf{G_\mathbf{\Psi}}$ and the identity matrix $\mathbf{I}_{N_c}$:

\begin{equation}\label{eq:L1:cc}
	\gamma(\mathbf{\Psi}) = \frac{1}{N} \min_{\mathbf{\Psi} \in R^{M \times N_c}} ||I_{N_c} - \mathbf{G_\mathbf{\Psi}}||^2_F
\end{equation}

where $||\cdot||_F$ denotes the Frobenius norm and $N := K \times (K - 1)$ is the total number of column pairs. Note that the optimization of only the average cross-correlation can result in large mutual coherence and is regularly prone to inaccurate recovery. In this context, Alemazkoor and Meidani (2018) \cite{Alemazkoor.2018} proposed a hybrid optimization criteria, which minimizes both the average cross-correlation $\gamma(\mathbf{\Psi})$ and the mutual coherence $\mu(\mathbf{\Psi})$:

\begin{equation}\label{key}
\argmin\left(f(\mathbf{\Psi})\right) = \argmin\left(\left(\frac{\mu_{i} -\min(\boldsymbol\mu)}{\max(\boldsymbol\mu) - \min(\boldsymbol\mu)}\right)^2 + \left(\frac{\gamma_i -\min(\boldsymbol\gamma)}{\max(\boldsymbol\gamma) - \min(\boldsymbol\gamma)}\right)^2\right)
\end{equation}

with $\boldsymbol\mu = (\mu_{1}, \mu_{2}, ..., \mu_{i})$ and $ \boldsymbol\gamma = (\gamma_1, \gamma_2, ..., \gamma_i)$


%
%

\subsection{Greedy algorithm to determine optimal sets of sampling points}				


		 
We used a greedy algorithm as shown in Algorithm \ref{alg:greedy} to determine L1-optimal sets of sampling points. In this algorithm, we generate a pool of $M_p$ samples and randomly pick an initial sample. In the next iteration, we successively add a sampling point and calculate the respective optimization criteria. After evaluating all possible candidates, we select the sampling point yielding the best criterion and append it to the existing set. This process is repeated until the sampling set has the desired size $M$.

\begin{algorithm}
\caption{Greedy algorithm to determine L1-optimal sets of sampling points}\label{alg:greedy}
\begin{algorithmic}[1]
	\State $\textit{create a random pool of } M_p \textit { samples}$
	\State $\textit{create the measurement matrix }\mathbf{\Psi_{pool}} \textit { of the samples}$
	\State $\textit{initiate } \mathbf{\Psi_{opt}} \textit{ as a random row } r \textit{ of } \mathbf{\Psi_{pool}}$
	\State $\textit{add row } r \textit{ to the added rows } r_{added} $ 
	\For{$i$ in (2, M)}
	\For{$j$ in $(1, M_p \textit{ without } r_{added})$}
	\State $\mathbf{\Psi_j} = \textit{row-concatenate } (\mathbf{\Psi_{opt}}, r_j)$
	\State $f_j = f(\mathbf{\Psi_j}) $
	
	\State $\textit{evaluate } f_j = f(\mathbf{\Psi_j})$
	\EndFor			
	\State $\textit{save } f_i = argmin(f_j) \textit{ for all j and } j_{best} $
	\State $\textit{add } r_{j_{best}} \textit{ to } r_{added}$
	\State $\mathbf{\Psi_{opt}} = \textit{row-concatenate } (\mathbf{\Psi_{opt}}, r_{j_{best}})$
	\EndFor
	\State$\textit{Return } \mathbf{\Psi_{opt}} \textit{ and } r_{added} = X_{best}$ 
\end{algorithmic}
\end{algorithm}

\section{Results}
\label{sec:Results}
The respective performances of the sampling schemes are thoroughly investigated based on four different scenarios. We compare the accuracy of the resulting GPCE approximation with respect to the original model and investigate the convergence properties and recoverability of the different sampling schemes. Following comparable studies \cite{Alemazkoor.2018}, we used uniformly distributed random variables in all examples and constructed GPCE bases using Legendre polynomials. This is the most general case, as any other input distribution can (in principle) be emulated by modifying the post-processing stage of the GPCE. The examples compare the sampling schemes on three theoretical test functions: (i) The Ishigami Function representing a low-dimensional problem that requires a high approximation order; (ii) the six-dimensional Rosenbrock Function representing a problem of medium dimension and approximation order; and (iii) the 30-dimensional Linear Paired Product (LPP) Function \cite{Alemazkoor.2018} using a low-order approximation. Finally, we consider a practical example, which consists of an electrode model used to measure the impedance of biological tissues for different frequencies. All sampling schemes were implemented in the open-source python package pygpc \cite{Weise.2020}, and the corresponding scripts to run the benchmarks are provided in the supplemental material. The sparse coefficient vectors were determined using the LARS-Lasso solver from scipy \cite{Virtanen.2020}.

A summary about the GPCE parameters for each test case is given in Table \ref{tab:Testfunctions:Overview}. For each test function, we successively increased the approximation order until an NRMSD of $\varepsilon<1^{-5}$ is reached. We assumed a very high number of sampling points. In this regard, we wanted to eliminate approximation order effects in order to focus on the convergence with respect to the number of sampling points.

For each sampling scheme and test case, we created a large set of sampling points that can be segmented in different sizes. For each case, we computed the associated GPCE approximation and calculated the normalized root mean square deviation (NRMSD) between the GPCE approximation $\tilde{y}$ and the solution of the original model $y$ using an independent test set containing $N_t = 10.000$ random sampling points. The NRMSD is given by:

\begin{equation}
\varepsilon=\frac{\sqrt{\frac{1}{N_t}\sum_{i=1}^{N_t}\left(y_i-\tilde{y}_i\right)^2}}{\max(\mathbf{y}) - \min(\mathbf{y})}
\end{equation}

We evaluated the average convergence of the NRMSD together with the success rate of each sampling scheme by considering $30$ repetitions. In addition, we quantified the convergence of the first two statistical moments, i.e. the mean and the standard deviation. The results are presented in the supplemental material. Reference values for the mean and standard deviation were obtained for each testfunction from $N=10^7$ evaluations from the original model functions.

\begin{table}[t]
	\caption{Overview of numerical examples.}
	\label{tab:Testfunctions:Overview}
	\renewcommand{\arraystretch}{1.3} 
	\centering
	\begin{tabular}{c c c c c c}
		\hline
		\textbf{Function} & \textbf{Problem}      & \textbf{Dim.} & \textbf{Order} & \textbf{Int. Order} & \textbf{Basisfunctions} \\ \hline
		    Ishigami      & Low dim., high order  & $2$           & $12$           & $2$                 & $91$                    \\
		   Rosenbrock     & Med. dim., med. order & $6$           & $5$            & $2$                 & $181$                   \\
		       LPP        & High dim., low order  & $30$          & $2$            & $2$                 & $496$                   \\
		    Electrode     & Application example   & $7$           & $5$            & $3$                 & $596$                   \\ \hline
	\end{tabular}
\end{table}

	
\subsection{Low-dimensional high-order problem (Ishigami function)}
As a first test case, we investigate the performance of the different sampling schemes considering the Ishigami function \cite{Ishigami.1990}. It is often used as an example for uncertainty and sensitivity analysis because it exhibits strong nonlinearity and nonmonotonicity. It is given by:

\begin{equation}\label{eq:Ishigami}
y = \sin(x_1) + a \sin(x_2)^2 + b x_3^4 \sin(x_1) 
\end{equation}

\begin{figure*}[tbh]
	\centering
	\includegraphics[width=0.49\textwidth]{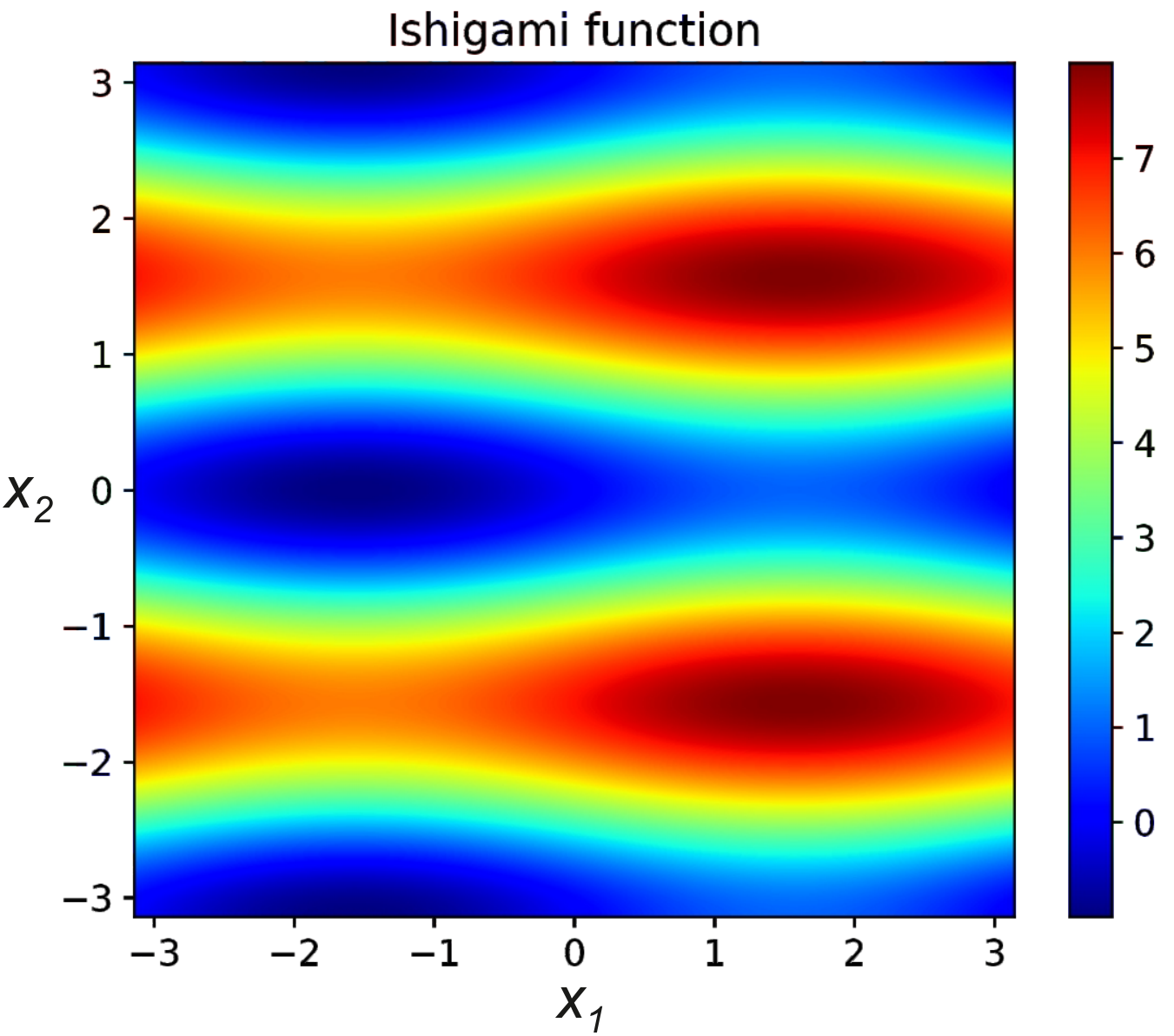}
	\caption[Ishigami]{Two-dimensional Ishigami function used to investigate the performance of different sampling schemes.}
	\label{fig:Testfunctions:Ishigami}
\end{figure*}

This example represents a low-dimensional problem requiring a high polynomial order to provide an accurate surrogate model. We defined $x_1$ and $x_2$ as uniformly distributed random variables and set $x_3=1$. The remaining constants are $a=7$ and $b=0.1$ according to \cite{Crestaux.2007} and \cite{Marrel.2009}. The approximation order was set to $p=12$, resulting in $N_c=91$ basis functions. We investigated the function in the interval $(-\pi, \pi)^2$ as shown in Fig. \ref{fig:Testfunctions:Ishigami}.

\begin{figure*}[!htbp]
	\centering
	\includegraphics[width=0.85\textwidth]{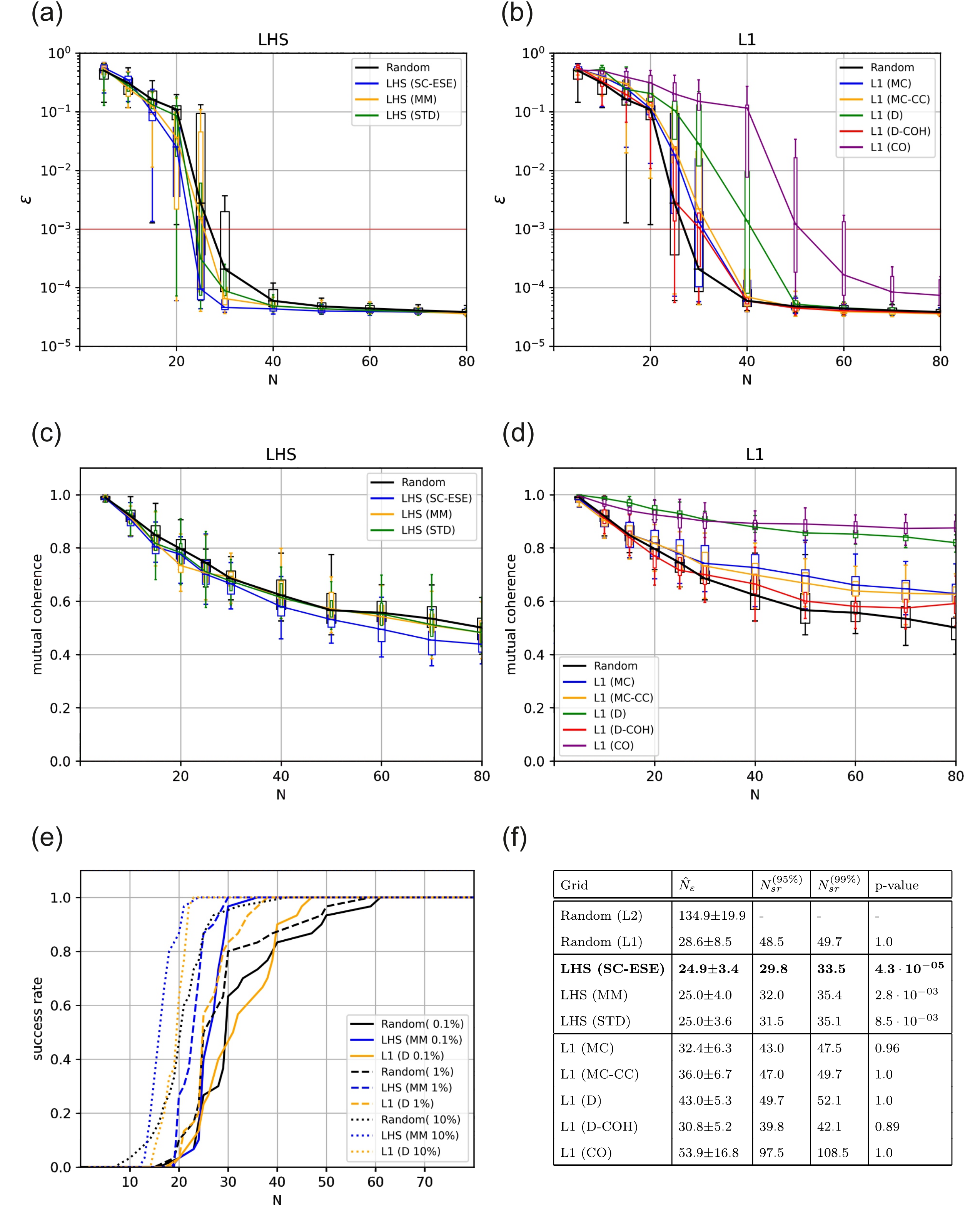}
	\caption{(a) and (b) Convergence of the NRMSD $\varepsilon$ with respect to the number of sampling points $N$ considering the Ishigami function. For reference, the convergence of the random sampling scheme is shown as a black line in each plot. (abbreviations: SC-ESE: stretched center enhanced stochastic evolutionary algorithm; MM: maximum-minimal distance; STD: standard; MC: mutual coherence; CC: cross-correlation; D: determinant optimal; D-COH: determinant-coherence optimal.); (c) and (d) mutual coherence of the gPC matrix; (e) success rate of the best converging grids for error thresholds of $0.1\%$, $1\%$, and $10\%$; (f) average number of sampling points needed to reach the error threshold of $0.1\%$, a success rate of $95\%$, $99\%$, and the associated p-value comparing if the grids require significantly lower number of sampling points than standard random sampling.}
	\label{fig:nrmsd:Ishigami}
\end{figure*}

The convergence results for the different sampling schemes are shown in Fig. \ref{fig:nrmsd:Ishigami}. It shows the dependence of the NRMSD $\varepsilon$ on the number of sampling points $N$ of the best sampling scheme from the LHS (Fig. \ref{fig:nrmsd:Ishigami}(a)) and L1-optimal sampling schemes (Fig. \ref{fig:nrmsd:Ishigami}(b)). The graphs of the error convergence consist of box-plots with whiskers of the error over different sampling sizes and are then connected with lines representing the median. Standard random sampling has the largest boxes and a black line on top for reference. We defined a target error level of $10^{-3}$, indicated by a horizontal red line, which corresponds to a relative error between the GPCE approximation and the original model function of 0.1\%. Additionally, the mutual coherence of the sampling sets is shown in Fig. \ref{fig:nrmsd:Ishigami}(c) and (d). The success rate of the best sampling schemes from both categories and standard random sampling is shown in Fig. \ref{fig:nrmsd:Ishigami}(e). In the Table shown in Fig. \ref{fig:nrmsd:Ishigami}(f), the median number of grid points required to reach that target error $\hat{N}_{\varepsilon}$ together with its standard deviation is shown. We also evaluated the median of the required number of sampling points for the random sampling scheme to determine the GPCE coefficients considering the L2 norm using the Moore-Penrose pseudo-inverse. All other evaluations have been performed using the LARS-Lasso solver. The success rates of the algorithms with the lowest 99\% recovery sampling size $\hat{N}_{\varepsilon}^{(99\% )}$ of each category of sampling schemes are marked in bold.


The sparsity of the model function greatly influences the reconstruction properties and hence the effectiveness of the sampling schemes. A GPCE approximation of the Ishigami function with an accuracy of $<10^{-5}$ requires $k=12$ out of the available $N_c=91$ coefficients ($13\%$).
Considering standard random sampling, the use of a conventional L2 solver requires $127$ sampling points to achieve a GPCE approximation with an error of less than $10^{-3}$ (see first row of Table in Fig. \ref{fig:nrmsd:Ishigami}(f)). In contrast, by using the L1 based LARS-Lasso solver, the number of required sampling points reduces to $31$, which serves as a baseline to compare the performance of the investigated sampling schemes.
By using the ESE enhanced LHS sampling scheme, the number of sampling points could be reduced to $25$, a substantial relative saving of $13$\% compared to standard random sampling. In the category of L1-optimal sampling, D-coherence optimal sampling schemes performed best; these schemes showed a slight increase in samples for the convergence. 

In addition to the average convergence of the sampling schemes, their reliability was calculated to evaluate their practical applicability. We have quantified reliability by calculating the number of sampling points required to achieve success rates of $95\%$ and $99\%$ in reaching the target error of $10^{-3}$; the reliability is determined by the relative number of repetitions required to reach the target error. Finally, we tested the hypothesis that the number of sampling points to reach the target error is significantly lower compared to standard random sampling. The Shapiro-Wilk-Test indicated that the error threshold distributions are not normally distributed. For this reason, we used the one-tailed Mann-Whitney U-test to compute the corresponding p-values. The generally good performance of LHS (SC-ESE) sampling is underpinned by a p-value of $4.3\cdot10^{-5}$. D-Coherence-optimal grids show a similar success rate and outperform standard random sampling (as measured by the numbers of sampling points required to achieve success rates of 95\% and 99\%) by factors of 9 and 8 respectively, signifying higher stability compared to standard random sampling on the Ishigami function.

Alongside the NRMSD, we calculated the mutual coherence of the GPCE matrix for each sampling scheme. It is shown in Fig. \ref{fig:nrmsd:Ishigami}(c) and (d). It can be seen that the mutual coherence is very stable around 0.7 for all LHS schemes considering a sampling size of 25. In contrast, L1-optimal grids show large variation. The coherence-optimal designs yield GPCE matrices with much higher coherences as defined in (\ref{eq:L1:mc}) compared to standard random sampling. D-Coherence-optimal sampling manages to reduce the mutual coherence the most after 25 samples; however, it shows very non-linear behavior for larger sampling sets, where it increases strongly above the level of random sampling.
	
\subsection{Medium-dimensional medium-order problem (Rosenbrock function)} \label{sec:Rosenbrock}
As a second test case, we used the $d$-dimensional generalized Rosenbrock function, also referred to as the Valley or Banana function \cite{Dixon.1978}. It is given by:

\begin{equation}\label{eq:Rosenbrock}
y = \sum_{i=1}^{d-1} 100\left(x_{i+1}-x_i^2\right)^2+\left(x_i-1\right)^2
\end{equation}

\begin{figure*}[tbh]
	\centering
	\includegraphics[width=0.49\textwidth]{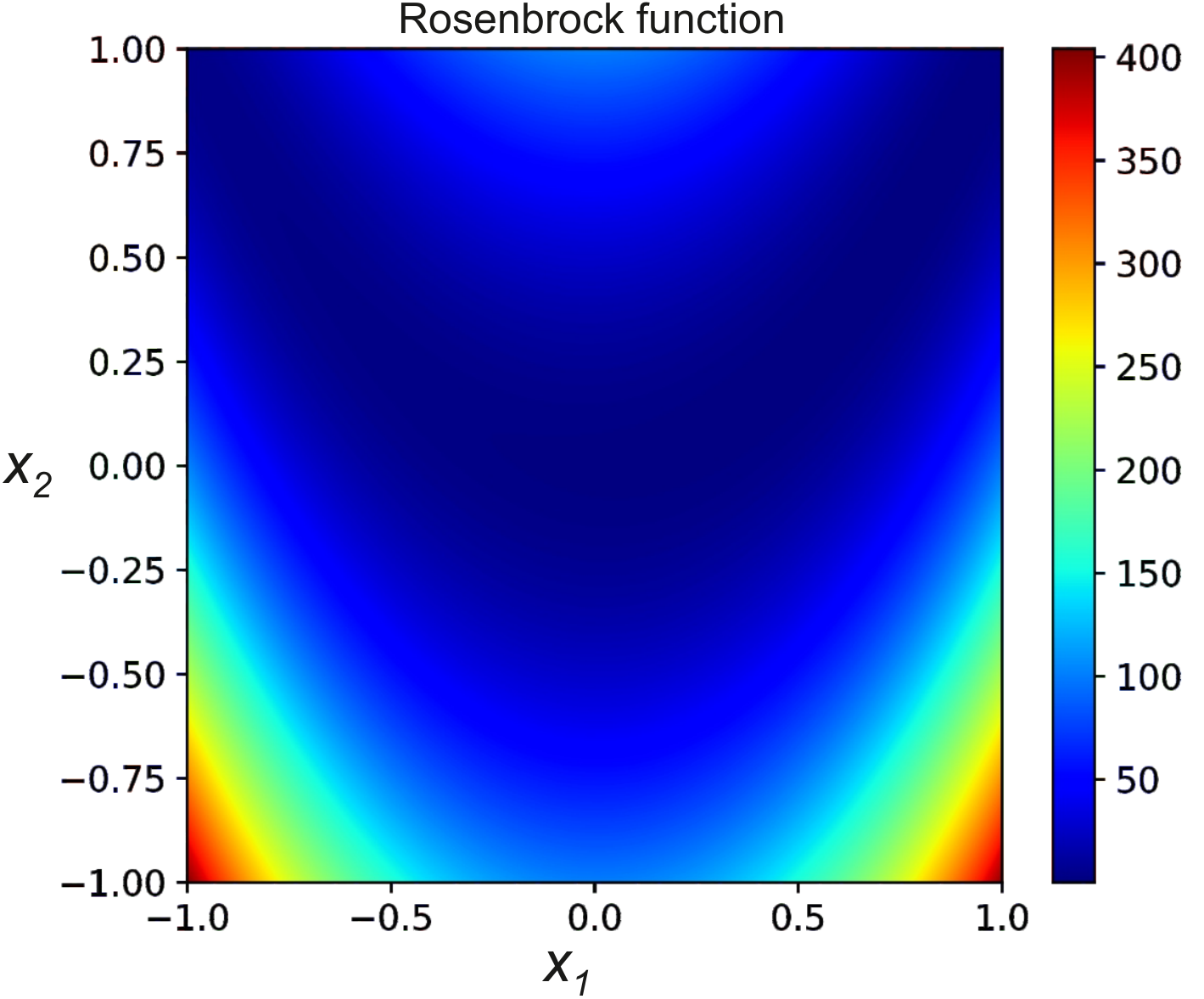}
	\caption[Rosenbrock]{Rosenbrock function in its two-dimensional form. In the present analysis, the Rosenbrock function of dimension $d=6$ is investigated.}
	\label{fig:Testfunctions:Rosenbrock}
\end{figure*}

The Rosenbrock function is a popular test problem for gradient-based optimization algorithms \cite{Molga.2015, Picheny.2013}. 
In Fig. \ref{fig:Testfunctions:Rosenbrock}, the function is shown in its two-dimensional form. The function is unimodal, and the global minimum lies in a narrow, parabolic valley. However, even though this valley is easy to approximate, it has more complex behavior close to the boundaries. To be consistent with our definitions of "low", "medium", and "high" dimensions and approximation orders in the current work, we classify this problem as medium-dimensional and requiring a moderate polynomial order approximation order to yield an accurate surrogate model. Accordingly, we defined the number of dimensions to be $d=6$ and set the approximation order to $p=5$ to ensure an approximation with an NRMSD of $\varepsilon<10^{-5}$ when using a high number of samples. This results in $N_c=181$ basis functions. The generalized Rosenbrock function is also used by Alemazkoor and Meidani \cite{Alemazkoor.2018} to compare the performance of different L1-optimal sampling strategies. We have used the same test function to make the results comparable and to be able to better integrate our study into the previous literature. 

\begin{figure*}[!htbp]
	\centering
	\includegraphics[width=0.85\textwidth]{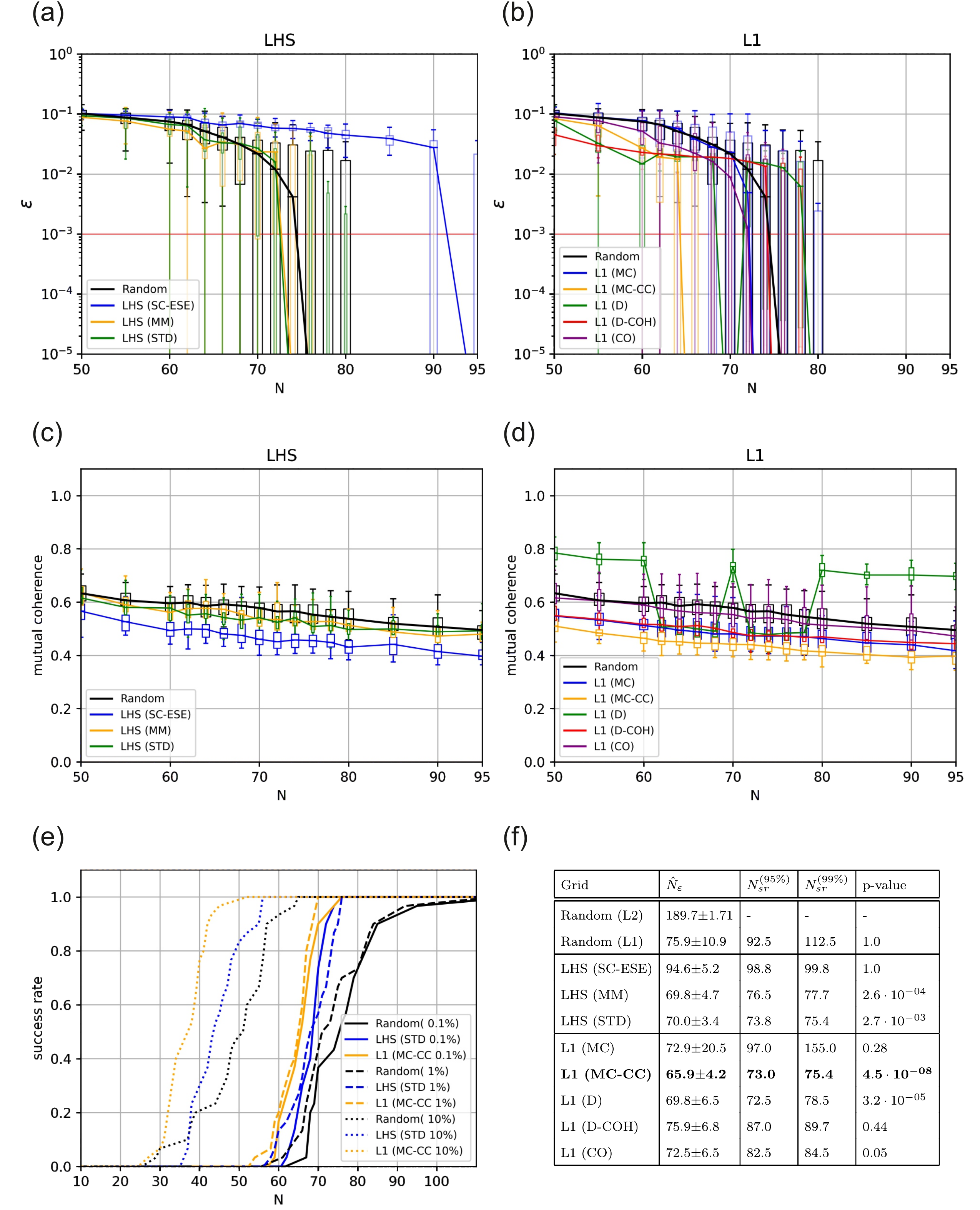}
	\caption{(a) and (b) Convergence of the NRMSD $\varepsilon$ with respect to the number of sampling points $N$ considering the Rosenbrock function. For reference, the convergence of the random sampling scheme is shown as a black line in each plot. (abbreviations: SC-ESE: stretched center enhanced stochastic evolutionary algorithm; MM: maximum-minimal distance; STD: standard; MC: mutual coherence; CC: cross-correlation; D: determinant optimal; D-COH: determinant-coherence optimal.); (c) and (d) mutual coherence of the gPC matrix; (e) success rate of the best converging grids for error thresholds of $0.1\%$, $1\%$, and $10\%$; (f) average number of sampling points needed to reach the error threshold of $0.1\%$, a success rate of $95\%$, $99\%$, and the associated p-value comparing if the grids require significantly lower number of sampling points than standard random sampling.}
	\label{fig:nrmsd:Rosenbrock}
\end{figure*}

The results of the error convergence for the investigated sampling schemes are shown in Fig. \ref{fig:nrmsd:Rosenbrock}(a) and (b). The mutual coherence of these algorithms is shown in Fig. \ref{fig:nrmsd:Rosenbrock}(c) and (d). The success rates of the best-performing sampling schemes from each category are shown in Fig. \ref{fig:nrmsd:Rosenbrock}(e), and the statistics are summarized in Table \ref{fig:nrmsd:Rosenbrock}(f). 
The Rosenbrock function can be exactly replicated by the polynomial basis functions of the GPCE using $k=23$ out of the $N_c=181$ available coefficients ($13\%$).
Random sampling in combination with L2-minimization requires $190$ sampling points to reach the target error of $10^{-3}$. In contrast, only $76$ samples are required when using the L1 based LARS-Lasso solver, which again serves as a baseline for comparison. The LHS (SC-ESE) algorithm is substantially less efficient than the other two LHS designs (STD and MM) in this test case. With LHS (MM) it is possible to achieve a reduction of sampling points by roughly 8\% compared to standard random sampling. From all investigated sampling schemes, MC-CC optimal grids performed best and required 13\% fewer sampling points than standard random grids. D-optimal sampling follows closely with a reduction of 8\%. However, for D-optimal, there is a very strong caveat connected to this measure which renders the value given by the table irrelevant. The median of the NRMSD for D-optimal sampling increases by orders of magnitude again after a sampling size of 75 is reached. It eventually drops below the error threshold again for sampling sizes close to 80 as seen in Fig. \ref{fig:nrmsd:Rosenbrock}(b), but this strong irregularity invalidates any statement regarding the significance of the error convergence for D-optimal sampling. This irregularity can also be observed in the mutual coherence as discussed in the following.

In terms of success rate, the sampling schemes differ considerably. Standard random sampling requires $N_{sr}^{95\%}=92.5$ and $N_{sr}^{99\%}=112.5$ sampling points to achieve success rates of $95\%$ and $99\%$ respectively. Standard LHS designs are more stable and require only $N_{sr}^{95\%}=84.5$ and $N_{sr}^{99\%}=84.9$ sampling points, respectively. D-optimal grids are very reliable and require only $N_{sr}^{95\%}=73.5$ and $N_{sr}^{99\%}=78.4$ sampling points. MC-CC grids, however, are able to surpass all other L1-optimal grids by achieving $N_{sr}^{99\%}=78.1$.

The mutual coherence of the measurement matrices for each algorithm is shown in Fig. \ref{fig:nrmsd:Rosenbrock}(c) and (d). It shows the same behavior for LHS sampling schemes as in case of the Ishigami function. This time, the sampling size of interest is larger, with about $80$ samples for the random sampling convergence. In this region, LHS (SC-ESE) shows the lowest mutual coherence at about $0.45$. Except for D-optimal sampling, the L1 sampling schemes are all able to reduce beyond the level of random sampling. This time, MC-CC sampling emerges as the leading design in that regard, which is then followed by MC and D-Coherence optimal designs. D-optimal designs show very irregular behavior connected to sampling sizes between $62$ to $68$ samples and $72$ to $78$, as can be seen in Fig. \ref{fig:nrmsd:Rosenbrock}(d). In both of these intervals, the mutual coherence drops briefly by about $0.3$ and then increases again to the initial level. Those two intervals also show tremendously lower NRMSD values then the surrounding sampling sizes in Fig. \ref{fig:nrmsd:Rosenbrock}(b).

\subsection{High-dimensional low-order problem (LPP function)}
As a third test case, we used the $d$-dimensional Linear Paired Product (LPP) function \cite{Alemazkoor.2018} assuming a linear combination between two consecutive dimensions:

\begin{equation}\label{eq:LPP}
y = \sum_{i=1}^{d}x_i x_{i+1}
\end{equation}

It has $d$ local minima except for the global one. It is continuous, convex and unimodal. In the present context, it is investigated having $d=30$ dimensions with an approximation order of $p=2$, resulting in $N_c=496$ basis functions. This test case represents high-dimensional problems requiring a low approximation order. This test function is also used by Alemazkoor and Meidani (2018) \cite{Alemazkoor.2018} but considering $d=20$ random variables.

\begin{figure*}[tbh]
	\centering
	\includegraphics[width=0.49\textwidth]{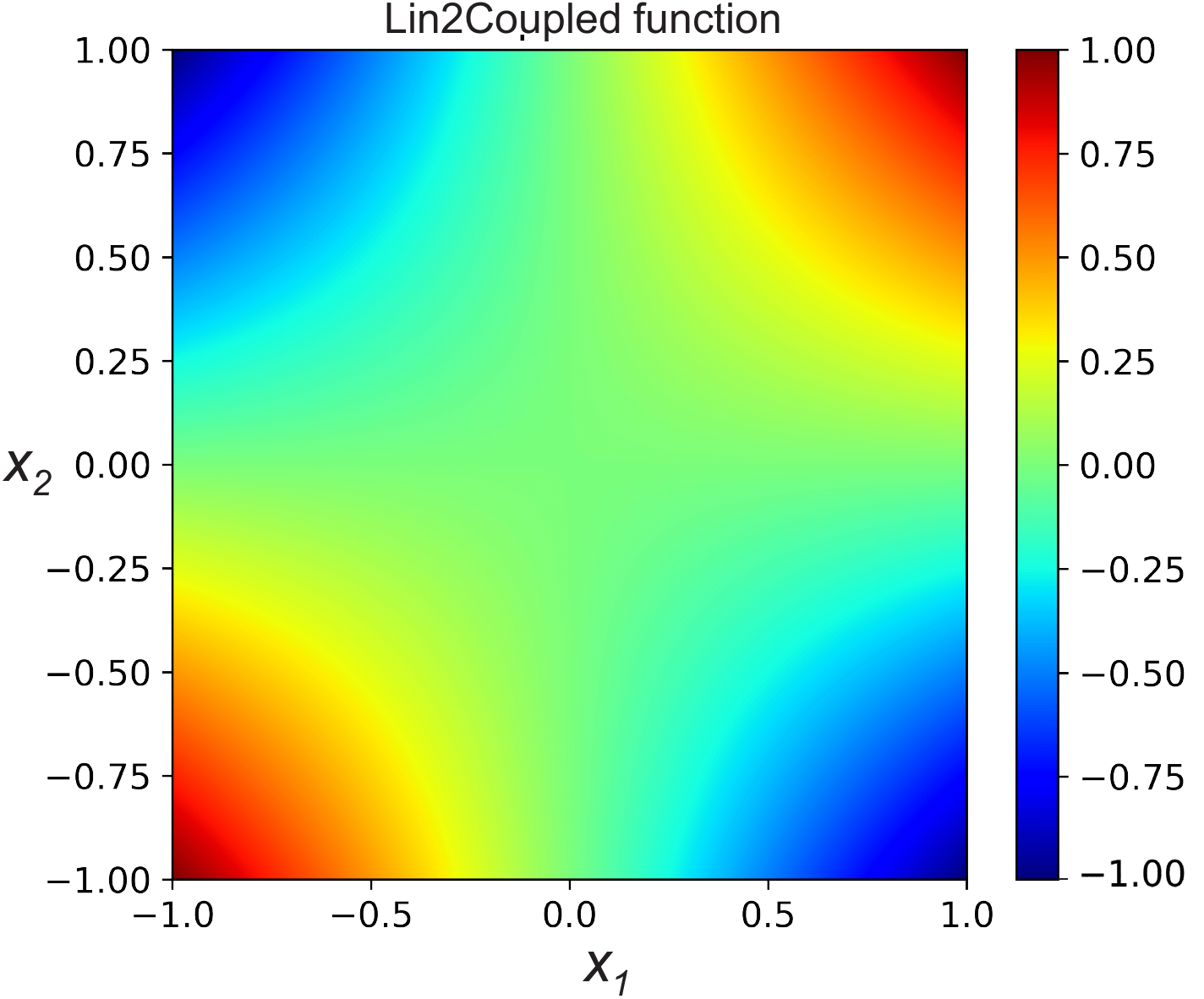}
	\caption[LPP]{Linear Paired Product (LPP) function in its two-dimensional form. In the present analysis, we investigated it with $d=30$ dimensions.}
	\label{fig:Testfunctions:Lin2Coupled}
\end{figure*}

\begin{figure*}[!htbp]
	\centering
	\includegraphics[width=0.85\textwidth]{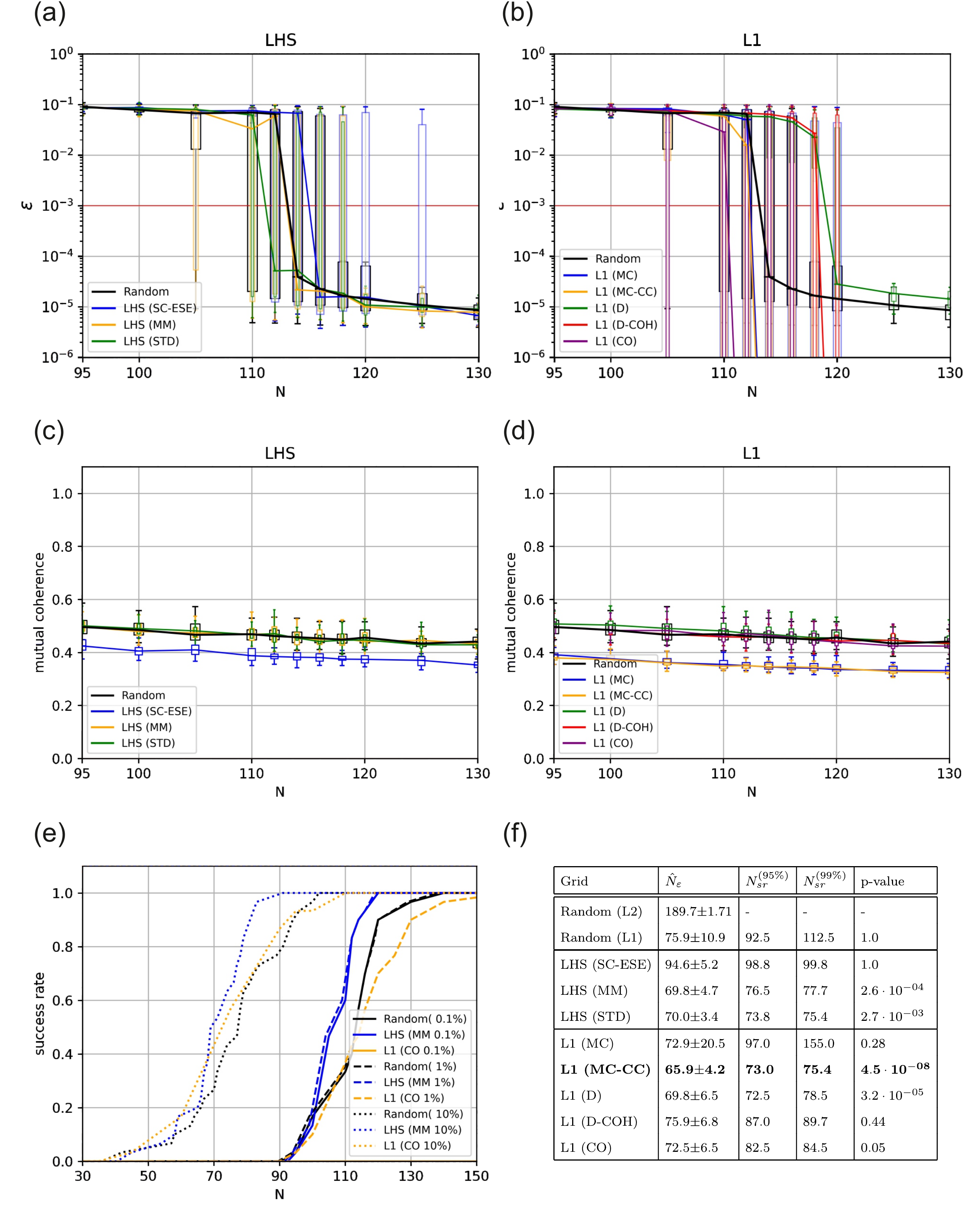}
	\caption{(a) and (b) Convergence of the NRMSD $\varepsilon$ with respect to the number of sampling points $N$ considering the Linear Paired Product (LPP) function. For reference, the convergence of the random sampling scheme is shown as a black line in each plot. (abbreviations: SC-ESE: stretched center enhanced stochastic evolutionary algorithm; MM: maximum-minimal distance; STD: standard; MC: mutual coherence; CC: cross-correlation; D: determinant optimal; D-COH: determinant-coherence optimal.); (c) and (d) mutual coherence of the gPC matrix; (e) success rate of the best converging grids for error thresholds of $0.1\%$, $1\%$, and $10\%$; (f) average number of sampling points needed to reach the error threshold of $0.1\%$, a success rate of $95\%$, $99\%$, and the associated p-value comparing if the grids require significantly lower number of sampling points than standard random sampling.}
	\label{fig:nrmsd:Lin2Coupled}
\end{figure*}

The error convergence for the different sampling schemes is shown in Fig. \ref{fig:nrmsd:Lin2Coupled}(a) and (b). The mutual coherence is visualized in \ref{fig:nrmsd:Lin2Coupled}(c) and (d), the success rates of the best-performing sampling schemes from each category are shown in Fig. \ref{fig:nrmsd:Lin2Coupled}(e), and the statistics are summarized in Table \ref{fig:nrmsd:Lin2Coupled}(d). 
This test function can be exactly replicated by the polynomial basis functions of the GPCE using $k=29$ out of $N_c=496$ available coefficients ($6\%$).
In this case, random sampling requires $190$ sampling points using L2-minimization and $110$ samples using L1-minimization. LHS designs showed similar convergence behavior to standard random sampling, with no reported improvement for LHS (MM) and LHS (STD) sampling and an increase in samples for LHS (SC-ESE). In contrast, L1-optimal designs did not manage to improve on the sampling count, with MC-CC and Coherence-Optimal sampling showing the best convergence rates. However, only D-optimal and D-Coherence optimal sampling increased by more then $4\%$, indicating very little variability between the sampling schemes. It can be observed that the variance in the range between $110$ and $120$ sampling points is very high for all sampling methods. The reason for this is that the LPP function is very sparse, and the convergence is mainly determined by the L1 solver. An additional sample point can lead to an abrupt reduction of the approximation error and a "perfect" recovery. This is often observed with L1 minimization.

To achieve the $95\%$ and $99\%$ success rates, standard random sampling requires $N_{sr}^{(95\%)}=122.5$ and $N_{sr}^{(99\%)}=128.5$ samples. The LHS (MM) algorithm performs slightly better and requires $N_{sr}^{(95\%)}=119.4$ and $N_{sr}^{(99\%)}=127$ samples. L1-optimal sampling schemes show tremendously weaker stability for this test-function. Here, only D-optimal designs reach the range of standard random sampling and the LHS variations, with $N_{sr}^{(95\%)}=128.1$ and $N_{sr}^{(99\%)}=129.6$.

The mutual coherence of the measurement matrices for this test case are shown in Fig. \ref{fig:nrmsd:Lin2Coupled}(c) and (d). The LHS (SC-ESE) shows the lowest coherence for the category of LHS grids, very similar to the previous example. L1 (MC) and L1 (MC-CC) designs display the lowest mutual coherence for L1-optimal designs, while the remaining L1-optimal sampling schemes are densely packed around the region slightly below random sampling. 

\subsection{Practical example (Probe impedance model)}
\begin{figure*}[tbh]
	\centering
	\includegraphics[width=0.9\textwidth]{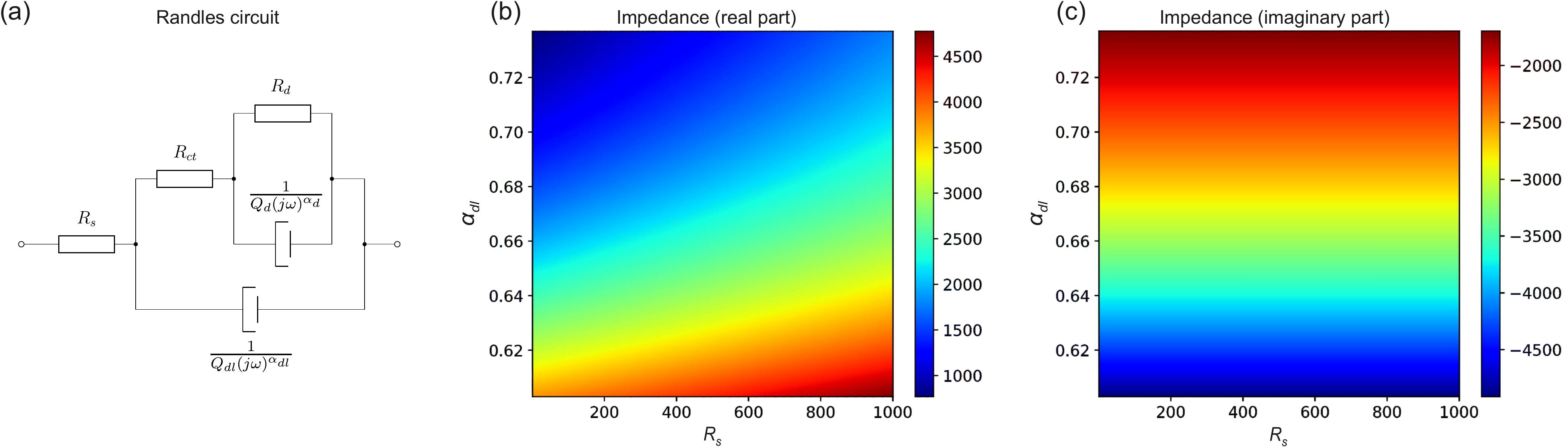}
	\caption[Electrode model]{Electrode impedance model: (a) Randles circuit; (b) Real part and (c) imaginary part of the electrode impedance as a function of $R_s$ and $\alpha_{dl}$. The remaining parameters are set to their respective mean values.}
	\label{fig:Testfunctions:Electrode}
\end{figure*}

The last test case is an application example from electrical engineering. The aim is to estimate the sensitivity of the intrinsic impedance of a probe used to measure the impedance of biological tissues for different frequencies. The model is shown in Fig. \ref{fig:Testfunctions:Electrode}(a) and consists of a Randles circuit that was modified according to the coaxial geometry of the electrode. The lumped parameters model the different contributions of the physical phenomena. The resistance $R_s$ models the contribution of the serial resistance of the electrolyte into which the electrode is dipped. The constant phase element $Q_{dl}$ models the distributed double layer capacitance of the electrode. The resistance $R_{ct}$ models the charge transfer resistance between the electrode and the electrolyte. The elements $Q_d$ and $R_d$ model the diffusion of charge carriers and other particles towards the electrode surface. The constant phase elements $Q_{dl}$ and $Q_d$ have impedances of $1/\left(Q_{dl}(j\omega)^{\alpha_{dl}}\right)$ and $1/\left(Q_{d}(j\omega)^{\alpha_{d}}\right)$, respectively. The electrode impedance, according to the Randles circuit shown in Fig. \ref{tab:Testfunctions:Electrode}(a), is given by:

\begin{equation}\label{eq:electrode}
\bar{Z}(\omega) = R_s + \left(Q_{dl}(j\omega)^{\alpha_{dl}} + \frac{1}{R_{ct}+\frac{R_d}{1+R_d Q_d (j\omega)^{\alpha_d}}}\right)^{-1}
\end{equation}

The impedance is complex-valued and depends on the angular frequency $\omega=2 \pi f$, which acts as an equivalent to the deterministic parameter $\mathbf{r}$ from eq. \ref{eq:ua:gPC_coeff_form}. A separate GPCE is constructed for each frequency. In this analysis, the angular frequency is varied between $1$ Hz and $1$ GHz with $1000$ logarithmically spaced points. The real part and the imaginary part are treated independently. The application example thus consists of $2000$ QOIs, for each of which a separate GPCE is performed. The approximation error is estimated by averaging the NRMSD over all QOIs. Accordingly, the impedance of the equivalent circuit depends on seven parameters, which will be treated as uncertain: $(R_s, R_{ct}, R_d, Q_d, \alpha_d, Q_{dl}, \alpha_{dl})$. They are modeled as uniformly-distributed random variables with a deviation of $\pm10$\% from their estimated mean values, with the exception of $R_s$, which was defined between $0 \Omega$ and $1 k\Omega$. The parameters were estimated by fitting the model to impedance measurements from a serial dilution experiment of KCl with different concentrations. The parameter limits are summarized in Table \ref{tab:Testfunctions:Electrode}. 
In preliminary investigations, we successively increased the approximation order until we reached an accurate surrogate model with an NRMSD of $\varepsilon<1^{-5}$. It was found that the parameters in this test problem strongly interact with each other, which explains the rather high order of approximation compared to the smooth progression of the real and imaginary parts in the cross sections shown in Fig. \ref{fig:Testfunctions:Electrode}. This means that (for example) when five parameters of first order interact with each other, the maximum GPCE order is reached, and this coefficient is significant compared to (for example) a fifth order approximation of a single parameter.

\begin{table}[t]
	\caption{Estimated mean values of the electrode impedance model, determined from calibration experiments, and limits of parameters.}
	\label{tab:Testfunctions:Electrode}
	\renewcommand{\arraystretch}{1.3} 
	\centering
	\begin{tabular}{c c c c}
		\hline
		\textbf{Parameter} & \textbf{Min.}   & \textbf{Mean}   & \textbf{Max.}   \\ \hline
		$R_s$        & $0$ $\Omega$    & $0.5$ k$\Omega$ & $1$ k$\Omega$   \\
		$R_{ct}$      & $9$ k$\Omega$   & $10$ k$\Omega$  & $1.1$ k$\Omega$ \\
		$R_d$        & $108$ k$\Omega$ & $120$ k$\Omega$ & $132$ k$\Omega$ \\
		$Q_{d}$       & $3.6 \cdot 10^{-10}$ F     & $4.0 \cdot 10^{-10}$ F     & $4.4 \cdot 10^{-10}$ F     \\
		$Q_{dl}$      & $5.4 \cdot 10^{-7}$ F      & $6 \cdot 10^{-7}$ F        & $6.6 \cdot 10^{-7}$ F      \\
		$\alpha_d$     & $0.855$         & $0.95$          & $1.0$           \\
		$\alpha_{dl}$    & $0.603$         & $0.67$          & $0.737$         \\ \hline
	\end{tabular}
\end{table}

\begin{figure*}[!htbp]
	\centering
	\includegraphics[width=0.85\textwidth]{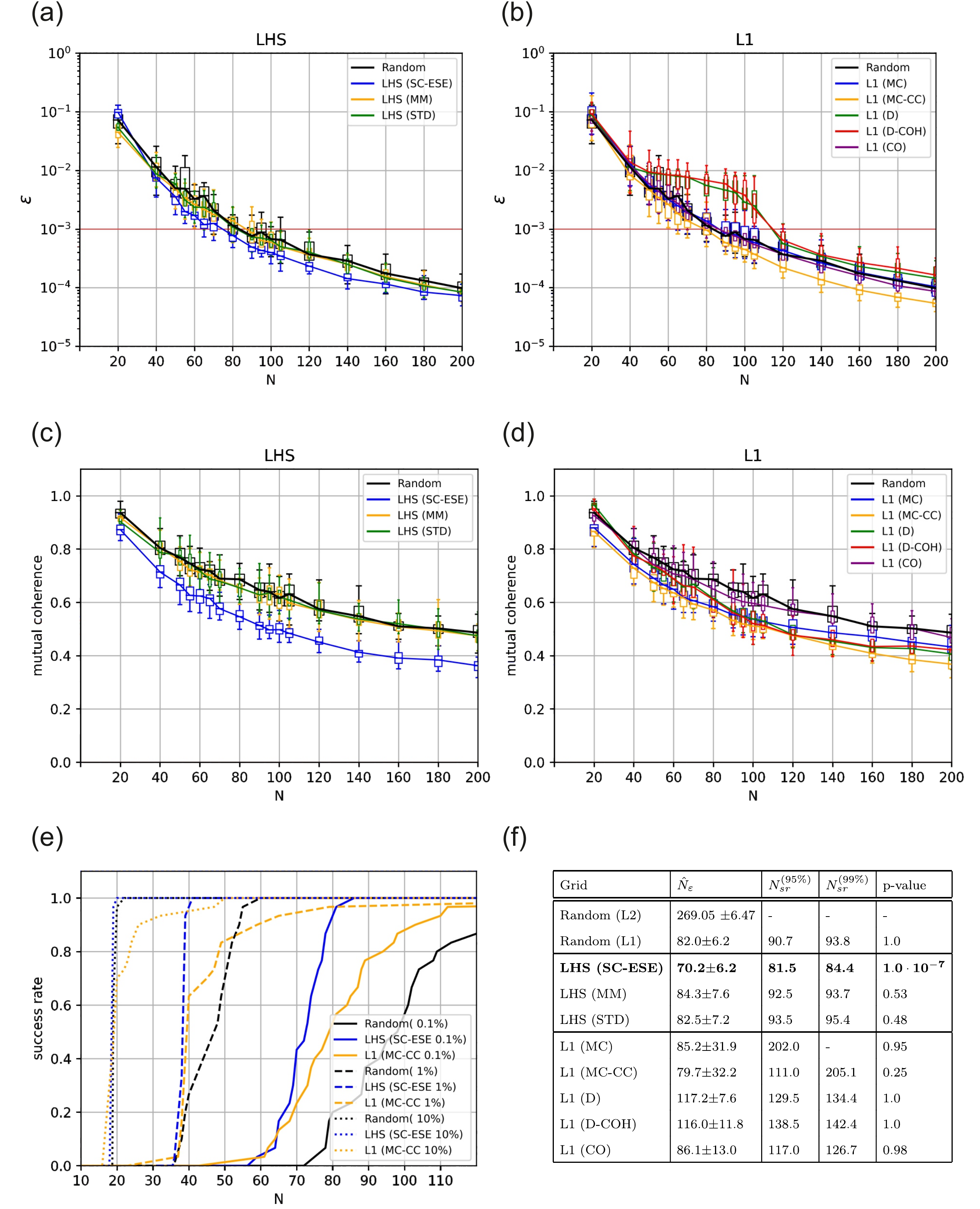}
	\caption{(a) and (b) Convergence of the NRMSD $\varepsilon$ with respect to the number of sampling points $N$ considering the probe impedance model. For reference, the convergence of the random sampling scheme is shown as a black line in each plot. (abbreviations: SC-ESE: stretched center enhanced stochastic evolutionary algorithm; MM: maximum-minimal distance; STD: standard; MC: mutual coherence; CC: cross-correlation; D: determinant optimal; D-COH: determinant-coherence optimal.); (c) and (d) mutual coherence of the gPC matrix; (e) success rate of the best converging grids for error thresholds of $0.1\%$, $1\%$, and $10\%$; (f) average number of sampling points needed to reach the error threshold of $0.1\%$, a success rate of $95\%$, $99\%$, and the associated p-value comparing if the grids require significantly lower number of sampling points than standard random sampling.}
	\label{fig:nrmsd:ElectrodeModel}
\end{figure*}

The results of the error convergence are shown in Fig. \ref{fig:nrmsd:ElectrodeModel}(a) and (b), and the corresponding mutual coherences are shown in \ref{fig:nrmsd:ElectrodeModel}(c) and (d). The success rates of the best-performing sampling schemes from each category are shown in Fig. \ref{fig:nrmsd:ElectrodeModel}(e), and the statistics are summarized in the Table Fig. \ref{fig:nrmsd:ElectrodeModel}(f).
The practical example consisting of the probe impedance model can be considered as non-sparse. It requires $k=500$ out of $N_c=596$ coefficients ($84\%$) to reach an accuracy of $<10^{-5}$.
Random sampling requires $269$ random samples to construct an accurate surrogate model using conventional L2 minimization. By using L1 minimization, the number of samples reduces to $82$. The LHS (SC-ESE) sampling scheme shows very good performance, requiring only $70.2$ samples to reach the target error, which corresponds to a decrease of $14$\%. In this test case, for L1-optimal sampling schemes, only MC-CC sampling managed to yield an improvement of the median by $2$ samples, yet they display a severe lack in stability as discussed in the next part.

Random grids require $N_{sr}^{(95\%)}=90.7$ and $N_{sr}^{(99\%)}=93.8$ samples to reach the desired success rates. Besides its good average convergence, LHS (SC-ESE) grids show significantly better stability, requiring only $N_{sr}^{(95\%)}=71.9$ and $N_{sr}^{(99\%)}=84.4$ samples, which corresponds to a decrease of $10$\% for both. A general lack of robust recovery is found for pure L1-optimal sampling schemes. Their success rates exceed those of random sampling, much like their median, rendering them inefficient on this test case. 

Fig. \ref{fig:nrmsd:ElectrodeModel}(c) and (d) shows the mutual coherences of the electrode impedance model. LHS (SC-ESE) still has the lowest mutual coherence in their category and may even show a lower mutual coherence than any L1-optimal design for single test runs. As seen previously, the lowest coherence for L1-optimal sampling can be observed in the case of MC-CC and MC sampling. Both of them form the bottom line of L1-optimal sampling in the region of $80$ samples; MC sampling rises above D and D-Coherence optimal sampling for larger sampling sizes.

\subsection{Average performance over all test problems}
The examined sampling schemes showed different strengths and weaknesses depending on the test problem. In order make general statements regarding their performance, we have weighted the error crossing $\hat{N}_{\varepsilon}$ and success rates $N_{sr}^{(95\%)}$ and $N_{sr}^{(99\%)}$ with the corresponding value of random sampling and averaged the results over all investigated test problems. The results are shown in Table \ref{tab:results:average}. It can be observed that LHS (SC-ESE) grids outperform random sampling in terms of average convergence and success rate for two of the test functions while maintaining a very competitive $N_{sr}^{(99\%)}$ on all test functions. They share this quality with the other LHS grids; however, the SC-ESE variation manages to score the largest sample reduction of $14.6$\% regarding the $\hat{N}_{\varepsilon}$ measure and about $35$\% on the two success rate measures. It is still closely trailed by the other two LHS sampling schemes. For the two higher-dimensional examples, however, LHS (STD) and LHS (MM) clearly show the most stable recovery success as seen in the $N_{sr}^{(99\%)}$ for all test functions, with a sample reduction of $15$\% on average. L1-optimal sampling schemes managed to achieve a significant sample reduction for the Rosenbrock function; specifically, MC-CC sampling is unrivalled regarding the $\hat{N}_{\varepsilon}$ decreases. However, the success rate measures are paralleled by LHS sampling. In terms of stability, only D-Coherence optimal sampling shows some perseverance, as its largest increase of samples for the $N_{sr}^{(99\%)}$ is $52$\% for the Electrode Model function. D-optimal sampling can be eliminated from this discussion since it showed a deeper-lying irregularity for the error convergence regarding the Rosenbrock function as discussed in Section \ref{sec:Rosenbrock}. Surprisingly, the remaining L1-optimal grids, namely MC, MC-CC and Coherence-optimal sampling, all exhibit sample increases of $100$\% and more for the $N_{sr}^{(99\%)}$ for one of the tested examples. For the test functions that were investigated in the context of a high-order GPCE approximation, no sampling scheme except MC-CC showed a reduction of samples, and all of them showed large increases in the amount of samples needed for the success rate targets. 
 
\begin{table}[!htbp]
	\caption{Relative and average number of grid points $\hat{N}_{\varepsilon}$ of different sampling schemes to reach an NRMSD of $10^{-3}$ with respect to standard random sampling using the LARS-Lasso solver (L1) considering all test functions. The columns for $N_{sr}^{(95\%)}$ and $N_{sr}^{(99\%)}$ show the number of samples needed to reach success rates of $95\%$ and $99\%$, respectively.}
	\label{tab:results:average}
	\setlength{\arrayrulewidth}{.2mm}
	\setlength{\tabcolsep}{4pt}		
	\renewcommand{\arraystretch}{1.85} 
	\centering
	\scriptsize
	\begin{tabular}{ |l||l|l|l|l|l|l|l|l|l|l|l|l|l|l|l| }
		\hline
		             & \multicolumn{3}{l}{Ishigami} \vline                             & \multicolumn{3}{l}{Rosenbrock} \vline                           & \multicolumn{3}{l}{LPP} \vline                                  & \multicolumn{3}{l}{Electrode} \vline                      & \multicolumn{3}{l}{Average (all testfunctions)}\vline                      \\ \hline
		Grid         & $\hat{N}_{\varepsilon}$ & $N_{sr}^{(95\%)}$ & $N_{sr}^{(99\%)}$ & $\hat{N}_{\varepsilon}$ & $N_{sr}^{(95\%)}$ & $N_{sr}^{(99\%)}$ & $\hat{N}_{\varepsilon}$ & $N_{sr}^{(95\%)}$ & $N_{sr}^{(99\%)}$ & $\hat{N}_{\varepsilon}$ & $N_{sr}^{(95\%)}$ & $N_{sr}^{(99\%)}$ & $\hat{N}_{\varepsilon}$ & $N_{sr}^{(95\%)}$ & $N_{sr}^{(99\%)}$ \\ \hline\hline
		Random (L1)  & 1                       & 1                 & 1                 & 1                       & 1                 & 1                 & 1                       & 1                 & 1                 & 1                       & 1                 & 1                 & 1                       & 1                 & 1                 \\ \hline\hline
		LHS (SC-ESE) & 0.869                   & 0.615             & 0.674             & 1.246                   & 1.068             & 0.887             & 1.018                   & 1                 & 1                 & 0.856                   & 0.898             & 0.900             & 0.997                   & 0.895             & 0.865             \\
		LHS (MM)     & 0.872                   & 0.660             & 0.712             & 0.919                   & 0.827             & 0.691             & 1                       & 0.972             & 0.966             & 1.027                   & 1.019             & 0.999             & 0.9545                  & 0.870             & 0.842             \\
		LHS (STD)    & 0.872                   & 0.649             & 0.706             & 0.921                   & 0.798             & 0.670             & 1                       & 0.987             & 0.983             & 1.006                   & 1.030             & 1.017             & 0.950                   & 0.866             & 0.844             \\ \hline\hline
		L1 (MC)      & 1.131                   & 0.887             & 0.956             & 0.961                   & 1.049             & 1.378             & 1.037                   & 1.282             & -                 & 1.039                   & 2.226             & -                 & 1.042                   & 1.361             & 1.167             \\
		L1 (MC-CC)   & 1.257                   & 0.969             & 1.000             & 0.868                   & 0.789             & 0.670             & 1.028                   & 1.235             & 1.449             & 0.972                   & 1.223             & 2.187             & 1.031                   & 1.054             & 1.327             \\
		L1 (D)       & 1.503                   & 1.026             & 1.048             & 0.92                    & 0.784             & 0.698             & 1.091                   & 1.079             & 1.082             & 1.429                   & 1.427             & 1.433             & 1.236                   & 1.079             & 1.065             \\
		L1 (D-COH)   & 1.076                   & 0.820             & 0.847             & 0.999                   & 0.941             & 0.797             & 1.082                   & 1.154             & 1.164             & 1.414                   & 1.526             & 1.518             & 1.143                   & 1.110             & 1.082             \\
		L1 (CO)      & 1.882                   & 2.010             & 2.183             & 0.955                   & 0.892             & 0.751             & 1.009                   & 1.132             & 1.160             & 1.05                    & 1.289             & 1.351             & 1.224                   & 1.331             & 1.361             \\ \hline
	\end{tabular}
\end{table}

\section{Discussion}\label{sec:Discussion}
We thoroughly investigated the convergence properties of space-filling sampling schemes, L1-optimal sampling schemes minimizing the mutual coherence of the GPCE matrix, hybrid versions of both, and optimal designs of experiment by considering different classes of problems. We compared their performance against standard random sampling and found great differences between the different sampling schemes. To the best of our knowledge, this study is currently the most comprehensive comparing different sampling schemes for GPCE.

Very consistently, a great reduction in the number of sampling points was observed for all test cases when using L1 minimization compared to least-squares approximations. The reason for this is that the GPCE basis in its traditional form of construction is almost always over-complete, and oftentimes not all basis functions and parameter interactions are required or appropriate for modeling the underlying transfer function. For this reason, the use of L1 minimization algorithms in the context of GPCE is strongly recommended as long as the approximation error is verified by an independent test set or leave-one-out cross validation.

The first three test cases can be considered "sparse" in the framework of GPCE. It is noted that these values are high in comparison to typical values in signal processing and do not fully meet the definition of sparse signals, where $||\mathbf{c}||_0 \ll N_c$. The sparse character is reflected in the shape of the convergence curves, which show a steep drop in the approximation error after reaching a certain number of sampling points. Non-sparse solutions (as in case of the probe impedance model) show a gradual exponential decrease of the approximation error.

The Ishigami function represents a class of problems where the QOI exhibits comparatively complex behavior within the parameter space. LHS designs, and especially LHS (SC-ESE), outperform all other investigated sampling schemes by taking advantage of their regular and space-filling properties, which ensures that model characteristics are covered over the whole sampling space. Similar benefits were observed for the probe impedance model.

D-optimal and D-coherence optimal designs were investigated in \cite{Hadigol.2018} by comparing their performance to standard random sampling using L2 minimization. In this context, the number of chosen sampling points had to be considerably larger than the number of basis functions, i.e. $N_g \gg N_c$. In the present analysis, we loosened this constraint and decreased the number of sampling points below the number of basis functions $N_g < N_c$. We were able to observe improved performance for the first three test cases (Ishigami function, Rosenbrock function, Linear Paired Product function), where the relative number of non-zero coefficients is between $6-13$\%. In the case of the non-sparse probe impedance model with a ratio of non-zero coefficients of $84\%$, D-optimal and D-coherence optimal designs were less efficient compared to standard random sampling.

In our analysis, we did not find a relationship between mutual coherence and error convergence. A good example is the excellent convergence of LHS algorithms and the comparatively high mutual coherence in the case of the Ishigami function (Fig. \ref{fig:nrmsd:Ishigami}) or the comparatively late convergence of mutual coherence optimized sampling schemes in the case of the Rosenbrock function (Fig. \ref{fig:nrmsd:Rosenbrock}). This is in accordance with the observations reported by Alemazkoor et al. (2018) \cite{Alemazkoor.2018}. However, it has been observed that minimizing maximum cross-correlation does not necessarily improve the recovery accuracy of compressive sampling. 
It can be concluded that both the properties of the transfer function and the sparsity of the model function greatly influence the reconstruction properties and hence the effectiveness of the sampling schemes.

We minimized the maximum cross-correlation of the measurement matrix $\mathbf{\Psi}$; however, this minimization demonstrated few benefits in reducing the number of sampling points to determine a GPCE approximation \cite{Li.1997, Elad.2007}.

All numerical examples were applied to uniformly distributed random inputs and Legendre polynomials, while in many real-world applications, different distributions may have to be assigned to each random variable. This requires the use of different polynomial basis functions, which would change the properties of the GPCE matrix. This can have a major influence on the performance of L1-optimal sampling schemes. In contrast, we expect fewer differences for LHS based grids because they only depend on the shape of the input pdfs and not additionally on the basis functions.

A large interest in the field of compressed sensing lies in identifying a lower bound on the sampling size needed for an accurate reconstruction. This depends on two factors: the properties of the measurement matrix and the sparsity of the solution. The formulation of general statements for GPCE matrices, which are constructed dynamically and depend on the number of random variables and type of input pdfs, is only possible with many preliminary considerations about their structure and type. Moreover, it requires additional work on the topic of sparsity estimation. Both topics are very recent and subject to current and future research. 

The sampling methods analyzed were based on the premise that the entire grid is created prior to the calculations. The importance of covering important aspects of the quantity of interest in the sampling space suggests the development of adaptive sampling methods. An iterative construction of the set of sampling points would benefit from information of already calculated function values. For this purpose, e.g. the gradients in the sampling points could be used to refine regions with high spatial frequencies. We believe that those methods would have great potential to further reduce the number of sampling points.


As a result, the convergence rate as well as the reliability could be increased considerably when using LHS or D-coherence optimal sampling schemes. Even though the LHS (SC-ESE) was enhanced in this paper to perform better in corner regions, it still performs less than optimally in cases where function features are close to the borders of the sampling region. In further investigations, it may become crucial to address the systematic bias in the optimization of the criterion used for the ESE algorithm by using the periodic distance as shown in section \ref{subsec:phi-limit} to fully remedy this shortcoming.

The advantages of the more advanced sampling methods over standard random sampling were even more pronounced when considering the first two statistical moments, i.e., the mean and standard deviation (see supplemental material).

We minimized the maximum cross-correlation of the measurement matrix $\mathbf{\Psi}$, but few benefits arose from reducing the number of sampling points to determine a GPCE approximation \cite{Li.1997, Elad.2007}. We could not observe that L1-optimal sampling schemes are superior to their competitors. It has also been repeatedly shown that $\mu$ may only be able to optimize the recovery in a worst case scenario and therefore acts as a lower bound on the expected error \cite{Elad.2007}. In this sense, we also could not observe a direct relationship between mutual coherence and error convergence. From our results, we conclude that the sampling points should be chosen to capture all properties of the model function under investigation rather than to optimize the properties of the GPCE matrix to yield an accurate surrogate model more efficiently. This is in contrast to the results reported by Alemazkoor et al. (2018) \cite{Alemazkoor.2018} in case of the LPP function. We could reproduce their results for mutual coherence optimal sampling. However, random sampling performed much better in our case than they reported. For this reason, they argued for an application of L1-optimized sampling schemes in case of high-dimensional sparse functions. They investigated the LPP function considering $20$ dimensions, whereas we considered $30$ random variables. Nevertheless, we could not reproduce their results for random sampling in this test case as well. A possible reason could be the use of a different L1 solver, which could lead to potential changes in effectiveness of certain algorithms and should therefore be taken into consideration when comparing the results.

\section*{Acknowledgments}
This work was supported by the German Science Foundation (DFG) (WE 59851/2); The NVIDIA Corporation (donation of one Titan Xp graphics card to KW). We acknowledge support for the publication costs by the Open Access Publication Fund of the Technische Universit\"at Ilmenau.

\bibliography{Manuscript_v3}	
	
\end{document}